\documentclass[a4paper,10pt]{article}
\usepackage{ifstyle}
\usepackage{graphicx}
\usepackage{epsfig} 

\usepackage[english]{babel}               
\usepackage[T1]{fontenc}
\usepackage[latin1]{inputenc}
\usepackage{amssymb}             
\usepackage{array}                 
\def\Xint#1{\mathchoice
{\XXint\displaystyle\textstyle{#1}}%
{\XXint\textstyle\scriptstyle{#1}}%
{\XXint\scriptstyle\scriptscriptstyle{#1}}%
{\XXint\scriptscriptstyle\scriptscriptstyle{#1}}%
\!\int}
\def\XXint#1#2#3{{\setbox0=\hbox{$#1{#2#3}{\int}$}
\vcenter{\hbox{$#2#3$}}\kern-.5\wd0}}

\def\dashint{\Xint-}
\newtheorem{theorem}{Theorem} 
\newtheorem{corol}{Corollary}

\newtheorem{lem}{Lemma}

\def\opn#1#2{\def#1{\operatorname{#2} } } 
\opn\Rm{Rm}
\opn\Ric{Ric} 
\opn\Rc{Rc}
\opn\Scal{Sc}
\opn\Tr{Tr}
\opn\Trac{Tr} 
\opn\dist{dist}
\opn\diam{Diam}
\opn\det{det} 
\opn\div{div}
\opn\Ker{Ker} 
\opn\exp{exp}
\opn\exph{exph}
\opn\Osc{Osc} 
\opn\Herm{Herm}
\opn\End{End} 
\opn\Hess{Hess} 
\opn\Vol{Vol}

\newcommand{\R}{\mathbb R}
\newcommand{\C}{\mathbb C}

\newcommand{\contract}{\mathrel{\kern-1.5pt\vrule width6.0pt height0.4pt depth0pt
                \vrule width0.4pt height4.0pt depth0pt}}
\newcommand{\retract}{\mathrel{\kern-1.5pt\vrule width0.4pt height4.0pt depth0pt
                          \vrule width6.0pt height0.4pt depth0pt}}

\newdimen\boxrulethickness \boxrulethickness=.07em
\newdimen\Openboxwidth \Openboxwidth=.72em
\newcommand{\Openbox}{%
 \leavevmode
 \hbox{%
   \hfil\vrule width\boxrulethickness
   \vbox to\Openboxwidth{%
     \advance\Openboxwidth -2\boxrulethickness
     \hrule height \boxrulethickness width\Openboxwidth\vfil
     \hrule height\boxrulethickness}%
   \vrule width\boxrulethickness\hfil
 }}
 			  
\begin{document}
%{\def\thefootnote{\relax}\footnote{\hskip-0.6cm{\bf{Mots-cl\'es}} : Connexions de Chern, Courbure de Chern, Vari\"ah\"ah presque complexes, Coordonn\"ahs presque complexes \\{\bf{Classification AMS}} : 32C35.}} 
\begin{center} 
\Large{\bf{Characterization of Einstein-Fano manifolds via the K\"ahler-Ricci flow}}
%The convergence of the K\"ahler-Ricci Flow over Einstein-Fano manifolds with trivial automorphism group}}
%\huge{\bf{Existence of K\"ahler-Einstein metrics over Fano contact manifolds with $b_2=1$}}
\\
\vspace{0.4cm}
\large{Nefton Pali}
\end{center} 
\begin{abstract}We explain a characterization of Einstein-Fano manifolds in terms of the lower bound of the density of the volume of the K\"ahler-Ricci Flow. This is a direct consequence of Perelman's uniform estimate for the K\"ahler-Ricci Flow and a $C^0$ estimate of Tian and Zhu.
%An application of our caracterisation is a free pluripotential Theory proof of the convergence of the K\"ahler-Ricci Flow over Einstein-Fano manifolds with trivial automorphism group. The more general case of the convergence of the K\"ahler-Ricci Flow over  solitonic Fano manifolds has been considered previously by Tian and Zhu.
\end{abstract}
%\tableofcontents 
\section{Introduction}
During his visit at the MIT in the spring of 2003 G.Perelman made the following surprising claim. Under the K\"ahler-Ricci Flow over a Fano manifold the normalized Ricci potential, its gradient and Laplacian, the diameter and the scalar curvature are uniformly bounded.
Perelman also gave a sketch of his proof. The proof uses in a crucial way the celebrated Perelman's no local collapsing result. The details have been filled out by Sesum and Tian in \cite{Se-Ti}. By using Perelman's result Tian and Zhu \cite{Ti-Zh} was able to prove the convergence of the K\"ahler-Ricci flow over solitonic Fano manifolds. In this way they partially prove the important Hamilton-Tian conjecture on the convergence of the K\"ahler-Ricci flow over Fano manifolds.
Perelman's spectacular result combined with the $C^0$ estimate of Tian and Zhu in \cite{Ti-Zh} implies directly the following characterization of Einstein-Fano manifolds in terms of the lower bound of the density of the volume of the K\"ahler-Ricci Flow.
\begin{theorem}
Let $X$ be a Fano manifold and $G$ be a compact maximal subgroup of the identity component of the group of automorphisms of $X$. Then $X$ admits a $G$-invariant K\"ahler-Einstein metric if and only if the K\"ahler-Ricci flow $(\omega_t)_t$ with $G$-invariant initial metric $\omega$ satisfies the uniform estimate $\omega_t^n\geq k\,\omega^n$, $k>0$ for all times $t\geq 0$.
\end{theorem}
This is an equivalent form of one of the main resulte in \cite{Ti-Zh}.
%We prove the required $C^0$ estimate without using any element of pluripotential Theory. This is possible by a carefull use of the expression of the K-energy along the K\"ahler-Ricci Flow.
In writing this fact we took also the occasion to give as much as possible an intrinsic flavor to the proof of the celebrated Yau's $C^2$ \cite{Yau} and Calabis's $C^3$-uniform estimates for the complex Monge-Amp\`ere equation in the case of the K\"ahler-Ricci flow (see also \cite{Cao}).
\\
The first step in proving Perelman's result consist in showing the boundedness of a normalizing constant which appears in the evolution formula of the Ricci potential. 
Perelman show this by using the monotonicity of his $\mu$ functional along the K\"ahler-Ricci flow. We realize that the boundedness  of this constant follows in a classical way by using the generalized Bochner-Kodaira Formula. This leads also to an intresting consequence. 
\begin{prop}\label{MonWkrfRP}%Perelman's ${\cal W}$ functional is increasing a
Along the K\"ahler-Ricci flow $\frac{d}{dt}\omega_t=\omega_t-\Ric_t=i\partial\bar\partial u_t$, 
\\
$\dashint_Xe^{-u_t}\omega^n_t=1$, Perelman's ${\cal W}$ functional of the Ricci potential $u_t$ with scale $\tau=1/2$ is increasing. Moreover the monotonicity is strict unless the flow is a soliton. 
\end{prop}
{\bf Acknowledgments.} The result presented in this note has been reported during the visit of the author in Princeton University. The author is very grateful to Professors Jean-Pierre Demailly, Joseph Kohn and Gang Tian who made possible his long visit in this institution. The author is especially grateful to Professor Gang Tian for all his continuous scientific support and encouragement.
%%%%%%%%%%%%%%%%%%%%%%%%%%%%%%%%%%%%%%%%%%%%%%%%%%%%%%%%%%%%%%%%%%%%%%%%%%%%%%%%%%%%%
%%%%%%%%%%%%%%%%%%%%%%%%%%%%%%%%%%%%%%%%%%%%%%%%%%%%%%%%%%%%%%%%%%%%%%%%%%%%%%%%%%%%%%%%%%%%%%%%%%%%%%%%%%%%%%%%%%%%%%%%%%%%%%%%%%%%%%%%%%%%%%%%%%%%%%%%%%%%%%%%%%%%%%%%APENDIX%%%%%%%%%%%%%%%%%%%%%%%%%%%%%%%%%%%%%%%%%%%%%%%%%%%%%%%%%%%%%%%%%%%%%%%%%%%%%%%%%%%%%%%%%%%%%%%%%%%%%%%%%%%%%%%%%%%%%%%%%%%%%%%%%%%%%%%%%%%%%%%%%%%%%%%%%%%%%%%%%%%%%%%%%%%%%%%%%%%%%%%%%%%
%\section{Background on complex differential geometry} 
\section{Curvature notions for the tangent bundle}
Let $(X,J)$ be a complex manifold of dimension $n$ equipped with a hermitian metric $\omega \in {\cal E}(\Lambda ^{1,1}_{_J}T_X^*)(X)$. We note by $D^{\omega}_{_J}=\partial^{\omega}+\bar{\partial}$ the Chern connection of the hermitian tangent bundle $(T_{X,J},h)$, where $h:=\omega (\cdot,J\cdot)-i\omega$ is the hermitian form on $T_{X,J}$ associated to $\omega$. We note by 
$$
{\cal C}_{\omega }(T_{_{X,J} }):=(D^{\omega}_{_J})^2
 \in {\cal E}(\Lambda ^{1,1}_{_{J}}T^*_{_X} \otimes_{_{\C} } \mbox{End}_{_{\C} }(T_{_{X,J}}))(X)
 $$
the Chern curvature form, which can also be given by the simpler formula 
\\
${\cal C}_{\omega }(T_{_{X,J} })\xi=\bar{\partial}\partial^{\omega}\xi$, for any germ of holomorphic vector field 
$\xi\in {\cal O}(T_{_{X,J} })_x$.
The Chern curvature
$
{\cal C}^{\omega}_{_{{X,J} } }  \in {\cal E}(\Herm ( T^{\otimes 2}_{_{X,J}}))(X)
$
is the hermitian form on the complex vector bundle $T^{\otimes 2}_{_{X,J}}$ defined by the formula 
\begin{eqnarray*} 
{\cal C}^{\omega}_{_{{X,J} } }(\xi _1\otimes\eta_1,\xi _2\otimes\eta_2)
&:=&
h({\cal C}_{\omega }(T_{_{X,J}})(\xi^{1,0} _1,\xi^{0,1} _2) \eta_1,\eta_2)
\\
\\
&=&-2i\omega({\cal C}_{\omega }(T_{_{X,J}})(\xi^{1,0} _1,\xi^{0,1} _2) \eta^{1,0} _1,\eta^{0,1}_2) 
\end{eqnarray*}
for all real vector fields $\xi _j,\eta_j\in{\cal E}(T_X)(U),\,j=1,2$ on some open subset $U$. The Griffiths curvature is defined by the formula
\begin{eqnarray*}
G^{\omega}_{_{{X,J}}}(\xi \otimes\eta):={\cal C}^{\omega}_{_{{X,J} } }(\xi \otimes\eta,\xi \otimes\eta)
=
\frac{1}{2} h(J{\cal C}_{\omega }(T_{_{X,J}})(\xi,J\xi) \eta,\eta).
\end{eqnarray*}
The fact that the Chern curvature is a hermitian form implies that the Griffiths curvature takes allways real values. Then we deduce the identity
\begin{eqnarray}\label{GriffIden}  
2G^{\omega}_{_{{X,J}}}(\xi \otimes\eta)=
\omega({\cal C}_{\omega }(T_{_{X,J}})(\xi,J\xi)\eta,\eta).
\end{eqnarray} 
If $\xi ,\eta \in {\cal O}(T_{X,J})(U)$ are holomorphic vector fields then the Griffiths curvature can be given by the simple formula
$$
G^{\omega}_{_{{X,J}}}(\xi \otimes\eta)=-\xi ^{1,0}.\,\xi ^{0,1}.\,|\eta|_{\omega}^2+|\partial^{\omega}_{{\xi} } \,\eta|_{\omega}^2\,.
$$
(see for example \cite{kob}).
Let $(z_1,...,z_n)$ be holomorphic coordinates and let 
\\
$(\zeta _k)_k\in {\cal O}(T^{1,0} _{_{X,J}} )^{\oplus n} (U)$ be a local holomorphic frame of the vector bundle  $T^{1,0} _{X,J}$. Consider the local expression of the metric 
$
\omega =\frac{i}{2} \sum_{ k,l}\omega _{k,\bar{l}}\,\zeta ^*_k\wedge\bar{\zeta}^*_l,
$
where the coefficients $\omega _{k,\bar{l}}$ satisfient the hermitian symmetry relation 
$
\overline{\omega _{k,\bar{l} }}=\omega _{l,\bar{k}}
$. 
We note by 
$
(\omega ^{k,\bar{l}}) =(\omega _{k,\bar{l}} )^{-1}
$
the inverse matrix of $(\omega _{k,\bar{l}})$, namely $\sum_t\omega ^{k,\bar{t}}\omega_{t,\bar{l}}=\delta_{k,l}$. If 
$\alpha\in\Lambda ^{p,q}_{_{J}}T^*_{_X} \otimes_{_{\C} } T^{1,0}_{_{X,J}}$ then we will note by 
$\alpha\otimes _{_{J}}\zeta^*_m:=\alpha\otimes\zeta^*_m+\overline{\alpha(\cdot)}\otimes\bar{\zeta}^*_m$.
With this notations the Chern curvature form is given locally by the expresion    
\begin{eqnarray*}
{\cal C}_{\omega }(T_{_{X,J} } )&=&\sum_{m =1}^n (\bar{\partial}\partial^{\omega}\zeta _m)\otimes _{_{J}}\zeta^*_m
=\sum_{l,m =1}^n C_{l,m}\otimes \zeta^*_m\otimes_{_{J}}\zeta _l
\\
&=&
\sum_{j,k,l,m =1}^n C^{j,\bar{k}}_{l,m}\,(dz_j\wedge d\bar{z }_k)\otimes \zeta^*_m\otimes_{_{J}}\zeta _l,
\end{eqnarray*}
with
\begin{eqnarray}\label{CoefCurv}  
C_{l,m}:=-\sum_{r=1}^n\,\Big(\partial\bar{\partial}\omega_{m,\bar{r}}-\sum_{s,t=1}^n\,\partial\omega_{m,\bar{s} }\wedge\omega ^{s,\bar{t}}\bar{\partial}\omega_{t,\bar{r}}\Big)\omega ^{r,\bar{l}}.  
\end{eqnarray}
%and 
%$C_{l,m}\otimes \zeta^*_m\otimes_{_{J}}\zeta _l:=C_{l,m}\otimes \zeta^*_m\otimes\zeta _l+
%\overline{C_{l,m}(\cdot)}\otimes \bar{\zeta }^*_m\otimes\bar{\zeta }_l$. 
The Chern curvature have the local expression
$$
{\cal C}^{\omega}_{_{{X,J} }}=\sum_{j,k,l,m=1}^n C_{j,l,\bar{k},\bar{m}}\,dz_j\otimes\zeta ^*_l\otimes d\bar{z}_k\otimes\bar{\zeta}^*_m, 
$$
where the coefficients $C_{j,l,\bar{k},\bar{m}}:=\sum_{h=1}^n C^{j,\bar{k}}_{h,l}\cdot \omega _{h,\bar{m}}$ satisfient the hermitian symmetry relation $\overline{C_{j,l,\bar{k},\bar{m}}}=C_{k,m,\bar{j},\bar{l}}$. The following lemma shows that the Chern curvature is the obstruction to the existence of holomorphic frames orthonormed at an order higher than one.
\begin{lem}\label{ChernSens} 
Let $(X,J)$ be a complex manifold of dimension $n$ equipped with a hermitian metric $\omega \in {\cal E}(\Lambda ^{1,1}_{_J}T_X^*)(X)$. Then for every point $x\in X$ and any $\omega(x)$-orthonormed frame $(e_k)_k\subset T^{1,0} _{_{X,J,x}}$
there exists holomorphic coordinates $(z_1,...,z_n)$ centered at $x$ and an holomorphic frame $(\zeta _k)_k\in {\cal O}(T^{1,0} _{_{X,J}} )^{\oplus n} (U_x)$, 
\\
$\zeta _k(x)=e_k$, in a neighborhood of $x$ such that the metric $\omega$ have the local expression
$$
\omega =\frac{i}{2} \sum_l\zeta ^*_l\wedge\bar{\zeta}^*_l
-\frac{i}{2} \sum_{ j,k,l,m} H_{l,\bar{m}}^{j,\bar{k} }\,z_j\bar{z}_k\,\zeta ^*_l\wedge\bar{\zeta}^*_m+O(|z|^3),
$$
where the coefficients $H_{l,\bar{m}}^{j,\bar{k} }$ satisfient the hermitian symmetry  $\overline{H_{l,\bar{m}}^{j,\bar{k}}}=H_{m,\bar{l}}^{k,\bar{j}}$.
Moreover for any such coordinates and frames the Chern curvatures have at the point $x$ the  expressions
\begin{eqnarray*}
&\displaystyle{
{\cal C}_{\omega }(T_{_{X,J} } )(x) =
\sum_{j,k,l,m =1}^n H^{j,\bar{k}}_{m,\bar{l}}\,(dz_j\wedge d\bar{z}_k)\otimes \zeta^*_m\otimes_{_{J}}\zeta _l,}&
\\
\\
&\displaystyle{
{\cal C}^{\omega}_{_{{X,J} }}(x)=\sum_{j,k,l,m=1}^n H^{j,\bar{k}}_{l,\bar{m}}\,dz_j\otimes\zeta ^*_l\otimes d\bar{z}_k\otimes\bar{\zeta}^*_m.}&
\end{eqnarray*}  
\end{lem} 
We define the Ricci tensor 
$\Ric_{_{J} }(\omega )\in {\cal E}(\Lambda ^{1,1}_{_J}T_X^*\cap \Lambda ^2_{_{\R} }T_X^* )(X)$ of the metric $\omega$ respect to the complex structure $J$ by the formula
$$
\Ric_{_{J} }(\omega ):=i\Trac_{_{\C} }  {\cal C}_{\omega }(T_{_{X,J}} )= i \,{\cal C}_{\omega }(K^{-1}_{_{X,J} } )\in 2\pi c_1(X)\,,
$$
where ${\cal C}_{\omega }(K^{-1}_{X,J})$ is the Chern curvature form of the anticanonical bundle $K^{-1}_{X,J}:=\Lambda^n_{_{\C} }T_{X,J}$. The scalar curvature $\Scal_{_{J} }(\omega)\in {\cal E}(X,\R)$ of $\omega$ respect to $J$ is defined by the formula
$$
\Scal_{_{J} }(\omega):=\Trac_{{\omega} }(\Ric_{_{J} }(\omega ))= \frac{2n\Ric_{_{J} }(\omega )\wedge \omega ^{n-1} }{\omega ^n}.  
$$ 
The fact that the Chern connection is invariant by scalar multiplications of the metric implies that $\Ric_{_{J} }(\lambda \omega )=\Ric_{_{J} }(\omega )$ for every real number $\lambda >0$. The Ricci curvature have the following local expression
$$
\Ric_{_{J} }(\omega )= i\sum_{1\leq j,k,l \leq n}C^{j,\bar{k} }_{l,l}\,\zeta ^*_j\wedge \bar{\zeta } ^*_k.
$$ 
We remind (cf. \cite{Dem}) that if $(L,h)\rightarrow (X,J)$ is a holomorphic hermitian line bundle and $\sigma \in {\cal O}(L\smallsetminus 0)(U)$ is a non vanishing holomorphic section over an open set $U$ then the local expression of the Chern curvature is given by the formula 
$$
{\cal C}_h(L)=-\partial \bar{\partial} \log\, |\sigma|^2_{h}
$$ 
on $U$.
If  $(\zeta _k)_k\in {\cal O}(T^{1,0} _{_{X,J}} )^{\oplus n} (U)$ is a local holomorphic frame of the  vector bundle  $T^{1,0} _{X,J}$ 
then $|\zeta _1\wedge...\wedge \zeta _n|^2_{\omega}=\det(\omega _{k,\bar{l}})$. We deduce that the local expression of the Ricci curvature is given by the formula
$$
\Ric_{_{J} }(\omega )=-i\partial \bar{\partial} \log\,\det(\omega _{k,\bar{l}}).
$$
If $\omega _1$ is an other $J$-invariant metric then we have the global identity
$$
\Ric_{_{J}}(\omega _1)-\Ric_{_{J}}(\omega )=-i\partial \bar{\partial} \log\left(\frac{\omega^n_1}{\omega ^n}\right). 
$$
%%%%%%%%%%%%%%%%%%%%%%%%%%%%%%%%%%%%%%%%%%%%%%%%%%%%%%%%%%%%%%%%%%%%%%%%%%%%%%%%%%%%%%%%%%%%%%%%%%%%%%%%%%%%%%%%%%%%%%%%%%%%%%%%%%%%%%%%%%%%%%%%%%%%%%%%%%%%%%%%%%%%%%%%%%%%%%%%%%%%%%%%%%%%%%%%%%%%%%%%%%%%%%%%%%%%%%%%%%%%%%%%%%%%%%%%%
\subsection{The K\"ahler case}
%%%%%%%%%%%%%%%%%%%%%%%%%%%%%%%%%%%%%%%%%%%%%%%%%%%%%%%%%%%%%%%%%%%%%%%%%%%%%%%
If $(X,J,\omega)$ is a K\"ahler manifold then the Chern connection coincides with the Levi-Civita connection of the $J$-invariant Riemannian metric $g\equiv g_{\omega ,J}$ associated to $\omega$. This implies that in the K\"ahler case the Chern curvature form coincides with the Riemann curvature form ${\cal R}_g $. In this case the Riemann curvature 
$$
R_g(\xi ,\eta,\mu ,\zeta )\equiv R_g(\xi \wedge\eta,\mu \wedge\zeta ):=g({\cal C}_{\omega }(T_{_{X,J} } )(\xi ,\eta)\zeta,\mu),
$$
($\xi ,\eta,\mu ,\zeta\in T_X$) is a smooth section of the vector bundle $S^2_{_{\R} }(\Lambda ^{1,1}_{_J}T_X^*\cap \Lambda ^2_{_{\R} }T_X^*)$. We consider also the $\C$-linear extension of the Riemann curvature on the complexified tangent bundle $T_X\otimes_{_{\R} }\C$. We have the equalities
\begin{eqnarray*} 
&\displaystyle{
R_g(\xi_1^{1,0},\xi ^{0,1}_2,\eta_1^{1,0},\eta^{0,1}_2)=-R_g(\xi_1^{1,0},\xi ^{0,1}_2,\eta^{0,1}_2,\eta^{1,0}_1)=}&
\\
\\
%&=&-g({\cal C}_{\omega }(T_{_{X,J} })(\xi_1^{1,0},\xi ^{0,1}_2)\eta^{1,0}_1,\eta^{0,1}_2)
%\\
%\\
&\displaystyle{
= i\omega ({\cal C}_{\omega }(T_{_{X,J} })(\xi_1^{1,0},\xi ^{0,1}_2)\eta^{1,0}_1,\eta^{0,1}_2)
=
-\frac{1}{2}{\cal C}^{\omega}_{_{{X,J} } }(\xi_1\otimes \eta_1, \xi_2 \otimes\eta_2)\,.}&
\end{eqnarray*}
The Riemann curvature have the following local expression in arbitrary holomorphic coordinates $(z_1,...,z_n)$
$$
R_g=\sum_{j,k,l,m=1}^nR_{j,\bar{k},l,\bar{m}}\,(dz_j\wedge d\bar{z}_k)\otimes(dz_l\wedge d\bar{z}_m ), 
$$
where the coefficients are given by the formula $2R_{j,\bar{k},l,\bar{m}}=-C_{j,l,\bar{k},\bar{m}}$
respect to the frame $(\zeta _j):=(\partial/\partial z_j)$. Using the formula \eqref{CoefCurv} respect to the frame $(\zeta _j)$ we deduce the expression
\begin{eqnarray} \label{CoefRiemCurv}
2R_{j,\bar{k},l,\bar{m}}=
\frac{\partial^2 \omega _{l,\bar{m}}}{\partial z_j\partial\bar{z}_k }  
- \sum_{s,t=1}^n\frac{\partial\omega_{l,\bar{s} }}{\partial z_j }\,\omega ^{s,\bar{t}}\, 
\frac{\partial\omega_{t,\bar{m}}}{\partial\bar{z}_k}.
\end{eqnarray} 
The facts that the Riemann curvature is real, is symmetric over $\Lambda ^{1,1}_{_J}T_X^*$ and the first Bianchi identity 
${\cal R}_g(\xi ,\eta)\mu +{\cal R}_g(\eta,\mu)\xi +{\cal R}_g(\mu,\xi)\eta=0$, are expressed in terms of the coefficients of the Riemann curvature by the symmetries
\begin{eqnarray*} 
&\displaystyle{
\overline{R_{j,\bar{k},l,\bar{m}}}=R_{k,\bar{j},m,\bar{l}} }& 
\\
&\displaystyle{
R_{j,\bar{k},l,\bar{m}} =R_{l,\bar{m},j,\bar{k}}  }&
\\
&\displaystyle{
R_{j,\bar{k},l,\bar{m}}  =R_{j,\bar{m},l,\bar{k}}},&
\end{eqnarray*}
(the second and last equality implies also 
$R_{j,\bar{k},l,\bar{m}} =R_{l,\bar{k},j,\bar{m}}$). By the other hand we see that the second and last equality follows immediately from the K\"ahler symmetries
$$
\frac{\partial \omega _{l,\bar{m} } }{\partial z_j }=\frac{\partial \omega _{j,\bar{m} } }{\partial z_l},\;\quad
\frac{\partial \omega _{l,\bar{m} } }{\partial \bar{z}_k }=\frac{\partial \omega _{l,\bar{k} } }{\partial \bar{z}_m}.
$$
\\
{\bf Holomorphic geodesic coordinates}.
In the K\"ahler case the conclusions of lemma \ref{ChernSens} holds for the frame $(\zeta _j):=(\partial/\partial z_j)$. In fact we have the following  strongest result.
\begin{lem}\label{KahlCurvGeod}
Let $(X,J,\omega)$ be a K\"ahler manifold of dimension $n$. Then for every point $x\in X$ and any $\omega(x)$-orthonormed frame $(e_k)_k\subset T^{1,0} _{_{X,J,x}}$ there exist holomorphic coordinates $(z_1,...,z_n)$ centered at $x$ such that 
$\frac{\partial}{\partial z_k}{\vphantom{z} }_{| _{_x}}=e_k$
and the metric $\omega$ have the local expression
$$
\omega =\frac{i}{2} \sum_l\,dz_l\wedge d\bar{z}_l
-\frac{i}{2} \sum_{ j,k,l,m} \,H_{l,\bar{m}}^{j,\bar{k} }\,z_j\bar{z}_k\,dz_l\wedge d\bar{z}_m+O(|z|^3),
$$
where the coefficients $H_{l,\bar{m}}^{j,\bar{k} }$ satisfient the symmetries  $\overline{H_{l,\bar{m}}^{j,\bar{k}}}=H_{m,\bar{l}}^{k,\bar{j}}$, and
$H_{l,\bar{m}}^{j,\bar{k}}=H_{j,\bar{m}}^{l,\bar{k}}=H_{l,\bar{k}}^{j,\bar{m}}$.
Moreover for any such coordinates the Chern curvatures have at the point $x$ the  expressions 
\begin{eqnarray*}
&\displaystyle{
{\cal C}_{\omega }(T_{_{X,J} } )(x) =
\sum_{j,k,l,m } H^{j,\bar{k}}_{m,\bar{l}}\,(dz_j\wedge d\bar{z}_k)\otimes dz_m\otimes_{_{J}}\frac{\partial }{\partial z_l},}&
\\
&\displaystyle{
{\cal C}^{\omega}_{_{{X,J} }}(x)=\sum_{j,k,l,m} H^{j,\bar{k}}_{l,\bar{m}}\,dz_j\otimes dz_l\otimes d\bar{z}_k\otimes d\bar{z}_m.  }&
\end{eqnarray*}  
\end{lem}  
The lemma shows that in the K\"ahler case the Chern curvature is the obstruction to the existence of holomorphic coordinates $(z_1,...,z_n)$ such that the frame $(\partial/\partial z_j)$ is orthonormed at an order higher than one. This coordinates are called geodesic holomorphic coordinates. %The reader can find the details of the proof in Demailly's book (\cite{Dem}, Chapter VI).
\\
It will also be usefull a more precise version of the lemma \ref{KahlCurvGeod}. We need first some notation. Consider the complex vector bundle 
$
F:=S^2_{_{\C} } \Lambda ^{1,1}_{_{J}}T^*_{_X} 
$
equipped with the connection $\nabla_F$ induced by the  complexified Levi-Civita connection. We have the following lemma.
\begin{lem}\label{KahlCurvGeod3}
Let $(X,J,\omega)$ be a K\"ahler manifold of dimension $n$. Then for every point $x\in X$ and any $\omega(x)$-orthonormed frame $(e_k)_k\subset T^{1,0} _{_{X,J,x}}$ there exist $\omega$-geodesic holomorphic coordinates $(z_1,...,z_n)$ centered at $x$ such that 
$\frac{\partial}{\partial z_k}{\vphantom{z} }_{| _{_x}}=e_k$
and the metric $\omega$ have the local expression $\omega =\frac{i}{2} \sum_{l,m}\omega_{l,\bar{m}}\,dz_l\wedge d\bar{z}_m$, with
$$
\omega_{l,\bar{m}} =\delta_{l,\bar{m}}
-\sum_{j,k} \,H_{l,\bar{m}}^{j,\bar{k} }\,z_j\bar{z}_k
-\sum_{p,j,k} \,\left(H_{l,\bar{m}}^{p,j,\bar{k} }\,z_pz_j\bar{z}_k+\overline{H_{m,\bar{l}}^{p,j,\bar{k}}}\,z_k\bar{z}_p\bar{z}_j\right)
+O(|z|^4),
$$
where the coefficients $H_{l,\bar{m}}^{j,\bar{k} }$ satisfient the symmetries of lemma \ref{KahlCurvGeod} and the coefficients $H_{l,\bar{m}}^{p,j,\bar{k}}$ are symmetric in the indexes $p,j,l$ and $k,m$.
Moreover for any such coordinates the Riemann curvature and its first covariant derivatives has at the point $x$ the  expressions 
\begin{eqnarray*}
&\displaystyle{
R_{\omega }(x) =\frac{1}{2}
\sum_{j,k,l,m } H^{j,\bar{k}}_{l,\bar{m}}\,(dz_j\wedge d\bar{z}_k)\otimes (dz_l\wedge d\bar{z}_m),}&
\\
&\displaystyle{
\nabla^{1,0}_F\,R_{\omega }(x) =\frac{1}{2}\sum_{j,k,l,m,p} H_{l,\bar{m}}^{p,j,\bar{k} }\,dz_p\otimes(dz_j\wedge d\bar{z}_k)\otimes (dz_l\wedge d\bar{z}_m), }&
\\
&\displaystyle{
\nabla^{0,1}_F\,R_{\omega }(x) =\frac{1}{2}\sum_{j,k,l,m,p}\overline{H_{m,\bar{l}}^{p,k,\bar{j}}}\,d\bar{z}_p\otimes(dz_j\wedge d\bar{z}_k)\otimes (dz_l\wedge d\bar{z}_m).}&
\end{eqnarray*}  
\end{lem}  
%$Proof$. The elementary proof is very similar to the proof of lemma \ref{KahlCurvGeod}. \\\\
{\bf The Bisectional curvature}. Consider now the bisectional curvature 
$$
b\sigma _g(\xi ,\eta):=R_g(\xi ,J\xi ,\eta,J\eta)=
4R_g(\xi^{1,0},\xi ^{0,1},\eta^{0,1},\eta^{1,0}),
$$ 
(the last equality follows from the identity 
$\xi \wedge J\xi=-2i\xi ^{1,0}\wedge \xi ^{0,1}$). We remark that the bisectional curvature coincides with the sectional curvature $\sigma_g(\xi ,\eta):=R_g(\xi ,\eta,\xi ,\eta)$ on complex lines, (in fact $\sigma_g(\xi ,J\xi)=b\sigma _g(\xi,\xi)$). The identity \eqref{GriffIden} shows that in the K\"ahler case the Griffiths curvature coincides (modulo a factor $2$) with the bisectional curvature. In the K\"ahler case the Riemann curvature is determined by the bisectional curvature. In fact the vector bundle $\Lambda ^{1,1}_{_J}T_X$ is generated over $\C$ by the vectors of type $\xi ^{1,0}\wedge \xi^{0,1}$, (see for example \cite{Dem}, Chapter III, sect 1). 
\\
%%%%%%%%%%%%%%%%%%%%%%%%%%%%%%%%%%%%%%%%%%%%%%%%%%%%%%%%%%%%%%%%%%%%%%%%%%%%%%%%%%%%%%%%%%%%%%%                     DEFINITION OF CURVATURE  OPERATOR
%%%%%%%%%%%%%%%%%%%%%%%%%%%%%%%%%%%%%%%%%%%%%%%%%%%%%%%%%%%%%%%%%%%%%%%%%%%%%%%%%%%%%%%%%%%%%%%%%%%%%%%%%%%%%%%%%%%%%%%%%%%%%%%%%%%%%%%%%%%%%%%%%%%%%%%%%%%%%%%%%%%%%%%%%%%%%%%%%%%%%                     DEFINITION OF CURVATURE  OPERATOR
%%%%%%%%%%%%%%%%%%%%%%%%%%%%%%%%%%%%%%%%%%%%%%%%%%%%%%%%%%%%%%%%%%%%%%%%%%%%%%%%%%%%%%%%%%
{\bf The Riemann curvature operator}. 
Let $G\in {\cal E}(S^2_{_{\R} }(\Lambda ^{1,1}_{_J}T_X^*\cap \Lambda ^2_{_{\R} }T_X^*))(X)$ be the induced metric over the real vector bundle 
$\Lambda ^{1,1}_{_J}T_X\cap \Lambda ^2_{_{\R} }T_X$. We will still note by $G\in {\cal E}(S^2_{_{\C} }(\Lambda ^{1,1}_{_J}T_X^*))(X)$ the $\C$-linear extension over the complexified vector bundle $\Lambda ^{1,1}_{_J}T_X$. Explicitly the metric $G$ is given by the formula
$$
G(u_1\wedge u_2,v_1\wedge v_2):=\det(g(u_k,v_l))_{k,l}=\omega (u_1,v_2)\cdot \omega (v_1,u_2),
$$
for any $u_1,v_1\in T^{1,0}_{_{X,J}}$ and  $u_2,v_2\in T^{0,1}_{_{X,J}}$. 
We remind now that the Riemann curvature operator 
$\Rm_g\in{\cal E}(\End_{_{\R} }(\Lambda ^{1,1}_{_J}T_X\cap \Lambda ^2_{_{\R} }T_X))(X)$ is defined by the formula
$$
G(\Rm_g(\xi \wedge \eta),\mu \wedge\zeta):=R_g(\xi \wedge \eta,\mu \wedge\zeta),
$$
for any $\xi ,\mu \in  T^{1,0}_{_{X,J}}$ and $\eta, \zeta \in T^{0,1}_{_{X,J}}$.
In local coordinates we find the expression
$$
\Rm_g=\sum_{j,k,s,t=1}^n \Rm_{j,\bar{k}}^{s,\bar{t}}
\Big(dz_j\wedge d\bar{z}_k\Big)\otimes \Big(\frac{\partial}{\partial z_s}  \wedge\frac{\partial}{\partial \bar{z}_t}\Big),  
$$
with
$$
\Rm_{j,\bar{k}}^{s,\bar{t}}=-4\sum_{l,m=1}^n\omega^{t,\bar{l} }  \omega^{m,\bar{s}} R_{j,\bar{k},l,\bar{m}}. 
$$
So in conclusion if put $H:=(\omega_{k,\bar{l}})$, we have the following synthetic expression$$
\Rm_g=2\sum_{s,t=1}^n 
H^{-1}\bar{\partial}(\partial H\cdot H^{-1})_{t,s}\otimes \Big(\frac{\partial}{\partial z_s}  \wedge\frac{\partial}{\partial \bar{z}_t}\Big)
$$
for the curvature operator.
The fact that $\Rm_g$ is a real operator implies the conditions 
$$
\overline{\Rm_{j,\bar{k}}^{s,\bar{t}}}=\Rm_{k,\bar{j}}^{t,\bar{s}}.
$$
%%%%%%%%%%%%%%%%%%%%%%%%%%%%%%%%%%%%%%%%%%%%%%%%%%%%%%%%%%%%%%%%%%%%%%%%%%%%%%%%%%%%%%%%%%%%%
%			    The Ricci Tensor %%%%%%%%%%%%%%%%%%%%%%%%%%%%%%%%%%%%%%%%%%%%%%%%%%%%%%%%%%%%%%%%%%%%%%%%%%%%%%%%%%%%%%%%%
{\bf The Ricci tensor}. We remind that in the Riemannian case the Ricci curvature $\Ric(g)\in {\cal E}(S^2_{_{\R} }T_X^*)(X)$ is defined by the formula 
$$
\Ric(g)(\xi ,\eta):=\Trac_{_{\R} } ({\cal R}_g(\cdot,\xi )\eta)
$$ 
for every $\xi ,\eta\in T_X$. If $(X,J,\omega)$ is a K\"ahler manifold and $g$ is the $J$-invariant Riemannian metric associated to $\omega$, then we have the formula 
$$
\Ric_{_{J} }(\omega )(\xi ,J\eta)=\Ric(g)(\xi ,\eta)
$$
for every $\xi ,\eta\in T_X$.  Let $(z_1,...,z_n)$ be $\omega$-geodesic coordinates centered in a point $x$ and set $\omega_0:=\frac{i}{2}\sum_k dz_k\wedge d\bar z_k$ and write $\Ric_{_{J} }(\omega )=i\sum_{k,l}R_{k\bar l}\,dz_k\wedge d\bar z_l$. Then we have the expansion
\begin{eqnarray}\label{geoMeanRic}
\omega^n=\left(1-\sum_{k,l}R_{k\bar l}(x)z_k\bar z_l\right)\omega^n_0+O(|z|^3)\,.
\end{eqnarray}
Starting from next section we will allways use Einstein's convention of sums.
%%%%%%%%%%%%%%%%%%%%%%%%%%%%%%%%%%%%%%%%%%%%%%%%%%%%%%%%%%%%%%%%%%%%%%%%%%%%%%%%%%%%%%%%%%%%%%%%%%%%%%%%%%%%%%%%%%%%%%%%%%%%%%
\section{The generalized Bochner-Kodaira formula for compact K\"ahler manifolds}
%%%%%%%%%%%%%%%%%%%%%%%%%%%%%%%%%%%%%%%%%%%%%%%%%%%%%%%%%%%%%%%%%%%%%%%%%%%%%%%%%%%%%%%%%%%%%%%%%%%%%%%%%%%%%%%%%%%
In writing this section we was inspired by \cite{Fu}. Let $(X,\omega )$ be a compact K\"ahler manifold of complex dimension $n$ and let
$$
\left<\alpha,\beta\right>_{\omega}:=\Tr_{\omega}(i\alpha\wedge \bar{\beta})/2=\frac{n\,i\alpha\wedge \bar{\beta}\wedge \omega^{n-1}}{\omega^n},
$$
be the induced hermitian product over the complex vector bundle $\Lambda ^{1,0}_{_{J}}T^*_{_X}$. Moreover let $h\in {\cal E}(X,\R)$ and $u\in {\cal E}(X,\C)$ be  smooth functions. Then the Laplacian 
$\Delta _{\omega,h }u:=\Delta _{\omega}u+2\left<\partial u,\partial h\right>_{\omega}$
is a self-adjoint differential operator respect to the inner product defined by the weighted volume form $e^h\omega^n$:
$$
(u,v)_{\omega,h }:=\int\limits_X u\bar{v}\,e^h\omega^n.
$$
In fact this follows from the identities $(\Delta _{\omega,h }u)e^h=-\Tr_{\omega}\left[i\bar{\partial}(e^h\partial u)\right]$
and
$$
-\int\limits_X i\bar{\partial}(e^h\partial u)\bar{v}\wedge \omega^{n-1}
=-\int\limits_X i\partial u\wedge \bar{\partial}\bar{v}\wedge e^h\omega^{n-1}
=\int\limits_X u\,i\partial(e^h\bar{\partial}\bar{v})\wedge \omega^{n-1}.
$$
We say that $\lambda\in \C$ is an eigenvalue of $\Delta _{\omega,h }$ if there exists a function $u\in {\cal E}(X,\C)$, not identically zero, such that $\Delta _{\omega,h}u+\lambda u=0$. Since
$$
-\int\limits_X(\Delta _{\omega,h }u)\,\bar{u}\,e^h\omega^n=\int\limits_X2|\partial u|^2_{\omega}\,e^h\omega^n,
$$ 
for any $u\in {\cal E}(X,\C)$, all the eigenvalues of $\Delta _{\omega,h }$ are nonnegative real numbers.
\begin{lem}{\bf (Generalized Bochner-Kodaira formula)}. Let $(X,\omega )$ be a compact K\"ahler manifold of complex dimension $n$ and let $u,\,h\in {\cal E}(X,\R)$ be smooth real functions. Then we have the Bochner type formula
\begin{eqnarray}\label{IntBochnh} 
\int\limits_X|\bar{\partial}\,\nabla^{1,0}_{\omega} u |^2_{\omega}\,e^h\omega ^n
&=&
-\int\limits_X\left<\partial\Delta _{\omega,h }u,\partial u\right>_{\omega}e^h\omega ^n\nonumber
\\
&-&
\int\limits_X\left(\Ric(\omega)-i\partial\bar{\partial}h\right)(\nabla_{\omega}u, J\nabla_{\omega}u)\,e^h\omega ^n. 
\end{eqnarray}     
\end{lem}  
$Proof$. Let $(z_1,...,z_n)$ be $\omega$-geodesic holomorphic coordinates with center a point $x$. By definition  of the $(2,0)$-component of the Hessian we have the identity
$
\nabla^{1,0}_{\omega}\partial u\,(\xi ,\eta)=\xi .\,\eta.\,u-(\nabla^{1,0}_{\omega,\xi }\,\eta ).\,u 
$
for every $(1,0)$-vector field $\xi ,\eta\in {\cal E}(T^{1,0}_{_{X,J}})(U)$ over an open set $U$.
By using the equality
$
\nabla^{1,0}_{\omega}\frac{\partial}{\partial z_l}=
\partial \omega _{l,\bar{j} }\,\omega ^{j,\bar{k}}\otimes\,\frac{\partial}{\partial z_k},      
$
we deduce the local expression
\begin{eqnarray*}
\nabla^{1,0}_{\omega}\partial u
&=&
\Big(u_{k,l}-\frac{\partial\omega _{l,\bar{j} } }{\partial z_k}\,\omega ^{j,\bar{r} }\,u_r \Big)\,dz_k\otimes dz_l 
\\
&=&
\left(u_{k,l}+C_{r,l}^{k,\bar{t}}\,\bar{z}_t\,u_r\right)\,dz_k\otimes dz_l+O(|z|^2).   
\end{eqnarray*}
Moreover the local expression $\nabla^{1,0} _{\omega}u=2u_{\bar{k} }\frac{\partial }{\partial z_k}+O(|z|^2)$ implies the local expression
\begin{eqnarray*}
\nabla^{1,0} _{\omega}u\contract\nabla^{1,0}_{\omega}\partial u
&=&
2(u_{k,l}u_{\bar{k} }+C_{r,l}^{k,\bar{t}}\,\bar{z}_t\,u_ru_{\bar{k} })\,dz_l+O(|z|^2)
\\
&=&
2(u_{k,l}u_{\bar{k} }+C_{t,l}^{k,\bar{r}}\,\bar{z}_t\,u_{\bar{k} }u_r)\,dz_l+O(|z|^2).
\end{eqnarray*}
We deduce the equality at the point $x$
\begin{eqnarray}\label{Boch-Rich} 
-\Tr_{\omega}\Big[ i\bar{\partial}(\nabla^{1,0} _{\omega}u\contract\nabla^{1,0}_{\omega}\partial u )\Big](x)
&=&
8(u_{k,l}u_{\bar{k} })_{\bar{l}} +8C_{l,l}^{k,\bar{r}}\,u_{\bar{k} }u_r\nonumber
\\
&=&
8(u_{k,l}u_{\bar{k} })_{\bar{l}} -2i\Ric(\omega)(\nabla^{1,0} _{\omega}u,\nabla^{0,1} _{\omega}u) \nonumber
\\
&=&
8(u_{k,l}u_{\bar{k} })_{\bar{l}}+\Ric(\omega)(\nabla_{\omega}u, J\nabla_{\omega}u)(x).\qquad
\end{eqnarray}
Consider now the trivial equalities at the point $x$
\begin{eqnarray*}
|\bar{\partial}\,\nabla^{1,0}_{\omega}u |^2_{\omega}(x)=8\,u_{k,l}u_{\bar{k},\bar{l}}
&=&8\,(u_{k,l}u_{\bar{k} })_{\bar{l}}-8\,u_{k,l,\bar{l} }\,u_{\bar{k}}
\\
&=&8(u_{k,l}u_{\bar{k} })_{\bar{l}}-\left<\partial\Delta _{\omega}u,\partial u\right>_{\omega}(x).
\end{eqnarray*}
Then using the equality \eqref{Boch-Rich} and the identity
\begin{eqnarray*}
\left<\partial\Delta _{\omega,h }u,\partial u\right>_{\omega}=\left<\partial\Delta _{\omega}u,\partial u\right>_{\omega}
+2\nabla^{1,0}_{\omega}\partial u(\nabla^{1,0}_{\omega}u,\nabla^{1,0}_{\omega}h)+i\partial\bar{\partial}h
(\nabla_{\omega}u, J\nabla_{\omega}u),
\end{eqnarray*}
we deduce the formula
\begin{eqnarray}\label{HorBoch}
|\bar{\partial}\,\nabla^{1,0}_{\omega}u |^2_{\omega}
&=&
-\left<\partial\Delta _{\omega,h }u,\partial u\right>_{\omega}
-\left(\Ric(\omega)-i\partial\bar{\partial}h \right)(\nabla_{\omega}u, J\nabla_{\omega}u)\nonumber
\\
&-&
\Tr_{\omega}\Big[i\bar{\partial}\left(\nabla^{1,0} _{\omega}u\,\contract\nabla^{1,0}_{\omega}\partial u \right)\Big]+2\nabla^{1,0}_{\omega}\partial u\left(\nabla^{1,0} _{\omega}u,\nabla^{1,0} _{\omega}h\right). \qquad
\end{eqnarray}
Moreover consider the equality
\begin{eqnarray*}
2n\,\bar{\partial}\Big[\left(i\nabla^{1,0} _{\omega}u\,\contract\nabla^{1,0}_{\omega}\partial u \right)\wedge e^h\omega^{n-1}\Big]
&=&
2n\,i\bar{\partial}\left(\nabla^{1,0} _{\omega}u\,\contract\nabla^{1,0}_{\omega}\partial u \right)\wedge e^h\omega^{n-1}
\\
&-&
2n\left(\nabla^{1,0} _{\omega}u\,\contract\nabla^{1,0}_{\omega}\partial u \right)\wedge i\bar{\partial}h\wedge e^h\omega^{n-1}.
\end{eqnarray*}
The last term is equal to
\begin{eqnarray*}
&-&2n\left(\nabla^{1,0} _{\omega}u\,\contract\nabla^{1,0}_{\omega}\partial u \right)\wedge \left(\nabla^{1,0}_{\omega}h\,\contract\omega\right)\wedge e^h\omega^{n-1}
\\
&=&-2\left(\nabla^{1,0} _{\omega}u\,\contract\nabla^{1,0}_{\omega}\partial u \right)\wedge  e^h\left(\nabla^{1,0}_{\omega}h\,\contract\omega^n\right)
\\
&=&-2\nabla^{1,0}_{\omega}\partial u\left(\nabla^{1,0} _{\omega}u,\nabla^{1,0} _{\omega}h\right)e^h\omega^n,
\end{eqnarray*}
so we have the formula
\begin{eqnarray*}
2n\,\bar{\partial}\Big[\left(i\nabla^{1,0} _{\omega}u\,\contract\nabla^{1,0}_{\omega}\partial u \right)\wedge e^h\omega^{n-1}\Big]
&=&
\Tr_{\omega}\Big[i\bar{\partial}\left(\nabla^{1,0} _{\omega}u\,\contract\nabla^{1,0}_{\omega}\partial u \right)\Big]e^h\omega^n
\\
&-&
2\nabla^{1,0}_{\omega}\partial u\left(\nabla^{1,0} _{\omega}u,\nabla^{1,0} _{\omega}h\right)e^h\omega^n.
\end{eqnarray*}
Then the formula \eqref{IntBochnh} follows from the formula \eqref{HorBoch} and the Stokes formula.\hfill$\Box$
\begin{corol}{\bf (Poincarr\'e type inequality)}.
Let $X$ be a Fano manifold of complex dimension $n$, let $\omega\in 2\pi c_1(X)$ be a K\"ahler metric and $h\in {\cal E}(X,\R)$ such that $\Ric(\omega)-\omega=i\partial\bar{\partial}h$. Set $V_h:=\int_Xe^h\omega^n$. Then for all smooth functions $\varphi\in {\cal E}(X,\R)$ we have the Poincarr\'e type inequality
\begin{eqnarray}\label{PoincIneq}
\int\limits_X|\partial \varphi|^2_{\omega}\,e^h\omega^n\geq
\int\limits_X\varphi^2e^h\omega^n-\frac{1}{V_h}\left(\;\int\limits_X\varphi\, e^h\omega^n\right)^2.
\end{eqnarray}
$Proof$. Let $u\in {\cal E}(X,\R)$ be an eigenfunction corresponding to the first eigenvalue  $\lambda _1>0$  of the Laplacian $\Delta _{\omega,h}$. Then the Bochner type formula \eqref{IntBochnh} implies the inequality
$$
\lambda _1\int\limits_X|\partial u|^2_{\omega}\,e^h\omega^n\geq 
\int\limits_X|\nabla_{\omega}u|^2_{\omega}\,e^h\omega^n=2\int\limits_X|\partial u|^2_{\omega}\,e^h\omega^n.
$$
The fact that $u$ can not be constant implies $\lambda _1\geq 2$. Consider now the function $\theta:=\varphi-\int_X\varphi\,e^h\omega^n/V_h$. Then the variational characterization of $\lambda _1$ implies the inequality
$$
\int\limits_X|\partial \theta|^2_{\omega}\,e^h\omega^n\geq \int\limits_X\theta^2e^h\omega^n,
$$
which implies the required Poincarr\'e type inequality \eqref{PoincIneq}.\hfill$\Box$
\end{corol}
%%%%%%%%%%%%%%%%%%%%%%%%%%%%%%%%%%%%%%%%%%%%%%%%%%%%%%%%%%%%%%%%%%%%%%%%%%%%%%%%%%%%%%%%%%%%%%%%%%%%%%%%%%%%%%%%%%%%%%%%%%%%%%%%%%%%%%%%%%%%%%%%%%%%%%%%%%%%%%%%%%%%%%%%%%%%%%%%%%%
\section{The K\"ahler-Ricci flow over Fano Manifolds}
%%%%%%%%%%%%%%%%%%%%%%%%%%%%%%%%%%%%%%%%%%%%%%%%%%%%%%%%%%%%%%%%%%%%%%%%%%%%%%%%%%%%%%%%%%%%%%%%%%%%%%%%%%%%%%%%%%%%%%%%%%%%%%%%%%%%%%%%%%%%%%%%%%%%%%%%%%%%%%%%%%%%%%%%%%%%%%%%%%%%%%%%%%%%%%%%%%%%%%%%%%%%%%%%%%%%%%%%%%%%%%%%%%%%%%%%%%%%%%%%%%%%%%%%%%%%%%%%%%%%
Let $X$ be a Fano manifold of complex dimension $n$ and let $\omega \in 2\pi c_1(X)$ be a K\"ahler metric. Let
$
{\cal P}_{\omega}:= \{\varphi \in {\cal E}(X,\R)\,|\,i\partial\bar{\partial}\varphi >-\omega\} 
$
be space of potentials
and define
$\omega _{\varphi}:=\omega +i\partial\bar{\partial}\varphi$ for every $\varphi \in {\cal P}_{\omega}$. 
The K\"ahler-Ricci flow is a family of K\"ahler metrics $(\omega _t)_t,$ solution of the evolution equation
\begin{eqnarray}\label{KRic-Flow}  
\frac{d}{dt}\omega _t=\omega _t-\Ric(\omega _t)
\end{eqnarray} 
with initial metric $\omega \in 2\pi c_1$. It was proved in \cite{Cao} that the K\"ahler-Ricci flow $(\omega _t)_t$ exists for all $t\in [0,+\infty)$ and $(\omega _t)_t\subset 2\pi c_1$. This is because solving the equation \eqref{KRic-Flow} is equivalent to solve the equation in terms of potentials 
\begin{eqnarray}\label{KRflowEqConst}
\dot{\varphi}_t=\log\,\frac{\omega _t ^n}{\omega ^n}+\varphi_t+c_t-h_{\omega},  
\end{eqnarray}
where $\varphi _t\in {\cal P}_{\omega},\,\omega _t=\omega +i\partial\bar{\partial}\varphi _t$, $\varphi_0=0$, $h_{\omega}\in {\cal E}(X,\R)$ is the the real smooth function defined by the conditions $\Ric(\omega)=\omega +i\partial\bar{\partial}h_{\omega}$, $\dashint_Xe^{h_{\omega}}\omega ^n=1$
and $c_t$ is a constant implying the normalization $\dashint_Xe^{-\dot{\varphi}_t}\omega^n_t=1$. We will allways consider the K\"ahler-Ricci flow equation with such normalization.
We remark that to find a solution $\varphi\in {\cal P}_{\omega}$ of the Einstein equation 
$
\Ric(\omega _{\varphi} )=\omega _{\varphi}, 
$
is equivalent to solve the equation 
$$
0=\log\,\frac{\omega _{\varphi}^n}{\omega ^n}+\varphi-h_{\omega}.
$$
This is also equivalent to the constant scalar curvature equation $\Scal(\omega _{\varphi})=2n$.
We prove now that the evolving metrics $\omega_t$ are $G$-invariant if the initial metric $\omega$ is $G$-invariant. Let $\Openbox_{\,t}:=\Delta _t-2\frac{\partial}{\partial t}$.
By deriving respect to a holomorphic vector field $\xi\in{\cal O}(T_X)(U)$ 
%or $\partial_{\bar k}:=\frac{\partial }{\partial \bar z_k}$, 
the K\"ahler-Ricci flow equation \eqref{KRflowEqConst} we find 
\begin{eqnarray}\label{C3Alfa}
\Openbox_{\,t}(\xi.\varphi_t)+2\xi .\varphi_t=(\Tr_{\omega}-\Tr_t)(L_{\xi}\,\omega)+2\xi . h_{\omega},
\end{eqnarray}
%where $\partial_k\omega:=\frac{i}{2}\partial_k\omega_{j\bar l}\,dz_j\wedge d\bar z_l$. 
This follows from the formula
$$
2\xi.\log\frac{\omega_t^n}{\omega^n}=\Tr_t L_{\xi}\,\omega_t-\Tr_{\omega}L_{\xi}\,\omega.
$$
Let prove this formula. Set $f_t:=\omega_t^n/\omega^n$. Then $L_{\xi}\,\omega_t^n=(\xi.f_t)\omega^n+f_t L_{\xi}\,\omega^n$.
%\begin{eqnarray*}
%L_{\xi}\,\omega_t^n&=&d(\xi \contract f_t\,\omega^n)=df_t\wedge (\xi \contract \omega^n)+f_td(\xi \contract \omega^n)\\
%&=&
%(\xi.f_t)\omega^n+nf_td(\xi \contract \omega)\wedge \omega^{n-1}.
%\end{eqnarray*}
So we get the equalities
\begin{eqnarray*}
nL_{\xi}\,\omega_t\wedge \omega^{n-1}_t=(\xi.f_t)\omega^n+nf_tL_{\xi}\,\omega\wedge \omega^{n-1},
\\
\xi.f_t=\frac{nL_{\xi}\,\omega_t\wedge \omega^{n-1}_t}{\omega^n}-f_t\frac{nL_{\xi}\,\omega\wedge \omega^{n-1}}{\omega^n},\quad
\\
2\,\xi.\log\frac{\omega_t^n}{\omega^n}=2\,\frac{\xi.f_t}{f_t}=\Tr_t L_{\xi}\,\omega_t-\Tr_{\omega}L_{\xi}\,\omega,\;\;
\end{eqnarray*}
which proves our formula. Let $\mathfrak{g}\subset H^0(T_{_X})$ be the (real) Lie algebra of $G$. We remark that a differential form $\alpha$ is $G$-invariant if and only if $L_{\xi}\,\alpha=0$ for all $\xi\in \frak{g}$.
Moreover the Ricci potential $h_{\omega}$ of any $G$-invariant metric $\omega$ is also $G$-invariant. So by applying \eqref{C3Alfa} with $\xi\in \frak{g}$ we find that the function $v_t:=\xi.\varphi_t$ is solution of the equation $\Openbox_{\,t}v_t=-2v_t$ with  initial data $v_0=0$. By uniqueness of the solutions we get $v_t=0$ and so the potential $\varphi_t$ is $G$-invariant for all times $t$.
\\
\\
{\bf K\"ahler-Ricci solitons.}
\\
Let $\omega$ be a K\"ahler metric and $u\in {\cal E}(X,\R)$ be a smooth real valued function. Then 
$\nabla_{\omega}u\contract \omega=-du\cdot J=-i\partial u+i\bar\partial u$ and 
$L_{\nabla_{\omega}u}\,\omega=d(\nabla_{\omega}u\contract \omega)=2i\partial\bar\partial u$. Let now $X$ be a Fano manifold and $\omega\in 2\pi c_1$ be a K\"ahler metric. Then 
$$
\omega-\Ric(\omega)=2i\partial\bar\partial u=L_{\nabla_{\omega}u}\,\omega.
$$
If $\nabla_{\omega}u\in {\cal O}(T_{_{X,J}})(X)$ then $\omega$ is called a K\"ahler-Ricci soliton.
We remind that $\nabla^{1,0}_{\omega}\partial u=0$ if and only if the vector field $\nabla_{\omega}u$ is holomorphic. So $\omega\in 2\pi c_1$ is a K\"ahler-Ricci soliton if and only if the Ricci potential $u\in {\cal E}(X,\R)$, $\omega-\Ric(\omega)=2i\partial\bar\partial u$ satisfies the equation 
$\nabla^{1,0}_{\omega}\partial u=0$. 
Let $(\Phi_t)_{t\in \R}$ be the 1-parameter Group of holomorphic automorphisms of $X$ induced by $\nabla_{\omega}u\in {\cal O}(T_{_{X}})(X)$. Let $\omega_t:=\Phi^*_t\omega =\omega +i\partial\bar{\partial}\varphi _t$ and $u_t=u\circ \Phi_t$. Then we get the K\"ahler-Ricci flow equation 
$$
\frac{d}{dt}\omega_t=\omega_t-\Ric(\omega_t)=2i\partial\bar\partial u_t,\quad\mbox{with}\quad \nabla^{1,0}_t\partial u_t=0.
$$
\\
{\bf Remark 1}. If the Futaki invariant $f_{2\pi c_1}$ is zero then all K\"ahler-Ricci solitons are K\"ahler-Einstein metrics.
In fact by definition of the Futaki invariant
$$
f_{2\pi c_1}(\nabla_{\omega} u)=-2\,\dashint\limits_X\nabla_{\omega} u\,.u\,\omega^n=-2\,\dashint\limits_X |\nabla_{\omega} u|^2_{\omega}\,\omega^n.
$$
Moreover this formula shows that the existence of a K\"ahler-Ricci solitons with $\nabla_{\omega} u\not =0$ implies the non existence of K\"ahler-Einstein metrics.
\\
\\
{\bf Remark 2}. Consider again a smooth real valued function $u\in {\cal E}(X,\R)$. Then 
$
L_{J\nabla_{\omega}u}\,\omega=d(J\nabla_{\omega}u\contract \omega)=0
$
since 
$
J\nabla_{\omega}u\contract \omega=-\omega(\nabla_{\omega}u, J\cdot)=-du.
$
Let now $X$ be a Fano manifold. For any K\"ahler metric $\omega$ and any smooth real vector field 
$\xi\in {\cal E}(T_{_X})(X)$ such that $L_{J\xi}\,\omega=0$ there exist a smooth real valued function 
$u\in {\cal E}(X,\R)$ such that $\xi=\nabla_{\omega}u$. In fact consider the decomposition $\xi=\xi'+\xi''$, with $\xi''=\overline{\xi'}\in {\cal E}(T^{0,1}_{_{X,J}})(X)$. Then
\begin{eqnarray*}
\;\;0=L_{J\xi}\,\omega=d(J\xi\contract \omega)=id(\xi'\contract \omega)-id(\xi''\contract \omega)
=
i\partial (\xi'\contract \omega)-i\bar\partial(\xi''\contract \omega),
\end{eqnarray*}
since $\bar\partial(\xi'\contract \omega)=0$ and $\partial(\xi''\contract \omega)=0$ by decomposition of the degree. We deduce also the equality $\partial (\xi'\contract \omega)=\bar\partial(\xi''\contract \omega)$.
The fact that $X$ is Fano and the equality $\bar\partial(\xi'\contract \omega)=0$ implient the existence of $u\in {\cal E}(X,\C)$ such that $\xi'\contract \omega=i\bar\partial u$ and by conjugation $\xi''\contract \omega=-i\partial\bar u$. Then the equality $\partial (\xi'\contract \omega)=\bar\partial(\xi''\contract \omega)$ implies $i\partial\bar\partial (u-\bar u)=0$, which means that the function $u$ can be chosen with real values.
\\
\\
{\bf Remark 3}. Let $(X,\omega)$ be a compact K\"ahler manifold such that $\omega-\Ric(\omega)=L_{\xi}\,\omega$, for some smooth real vector field $\xi\in {\cal E}(T_{_X})(X)$. So $L_{\xi}\,\omega=d(\xi\contract \omega)$ is a real $d$-exact $(1,1)$-form. By Hodge Theory there exist $u\in {\cal E}(X,\R)$ such that $L_{\xi}\,\omega=i\partial\bar\partial u$. So we deduce that $X$ is a Fano manifold and $\omega\in 2\pi c_1$. Using remark 2 we find that a K\"ahler metric $\omega$ over a compact K\"ahler manifold $X$ is a K\"ahler-Ricci soliton if and only if there exist a real holomorphic vector field $\xi\in {\cal O}(T_{_{X}})(X)$ such that $\omega-\Ric(\omega)=L_{\xi}\,\omega$ and $L_{J\xi}\,\omega=0$.
Moreover the holomorphic vector field $\xi$ is uniquely determined by the metric $\omega$, since it is uniquely determined by the Ricci potential of $\omega$.
\\
\\
{\bf Perelman's uniform estimates for the K\"ahler-Ricci flow.}
\\
%%%%%%%%%%%%%%%%%%%%%%%%%%%%%%%%%%%%%%%%%%%%%%%%%%%%%%%%%%%%%%%%%%%%%%%%%%%%%%%%%%%%%%%%%%%%%%%%%%%%%%%%%%%%%%%%%%%%%%%%%%%%%%%%%%%%%%%%%%%%%%%%%%%%%%%%%%%%%%%%%%%%%%%%%%%%%%%%%%%%%%
We have the following fundamental result due to Perelman.
\begin{theorem}{\bf (Perelman)}
Over a Fano manifold $X$ of complex dimension $n$, the K\"ahler-Ricci flow $\frac{d}{dt}\omega_t=\omega_t-\Ric_t=i\partial\bar\partial \dot{\varphi}_t$  satisfies the uniform estimates $|\dot{\varphi}_t|,\,|\nabla_t\,\dot{\varphi}_t|_t,\,|\Delta_t\dot{\varphi}_t|$, 
$
\diam_t(X),\,\Scal_t \leq C$, where the Ricci potential $\dot{\varphi}_t$ is normalized by the condition 
$\dashint_Xe^{-\dot{\varphi}_t}\omega^n_t=1$.
\end{theorem}
Set $u_t:=\dot{\varphi}_t,\;a_t:=\dot{c}_t$. Then time deriving the K\"ahler-Ricci flow equation \eqref{KRflowEqConst}, we find the identity 
\begin{eqnarray}\label{EvolUfunc}
\Openbox_{\,t}u_t=-2u_t-2a_t.
\end{eqnarray}
The first step in proving Perelman's theorem consist in showing the boundedness of the constant $a_t$. We give a proof here.
 %\cite{Per} 
By deriving the integral normalization $\dashint_Xe^{-u_t}\omega^n_t=1$ and using the equation \eqref{EvolUfunc}, we find the equalities
\begin{eqnarray*}
0=\frac{d}{dt}\,\dashint\limits_Xe^{-u_t}\omega^n_t
&=&
-\dashint\limits_X\dot{u}_t\,e^{-u_t}\omega^n_t
+2^{-1}\dashint\limits_X\Delta_t u_t\,e^{-u_t}\omega^n_t
\\
&=&
-a_t\dashint\limits_Xe^{-u_t}\omega^n_t-\dashint\limits_X u_t\,e^{-u_t}\omega^n_t,
\end{eqnarray*}
which implient
$$
a_t=-\dashint\limits_X u_t\,e^{-u_t}\omega^n_t\geq-\dashint\limits_X (u_t)_+\,e^{-u_t}\omega^n_t,
$$
where $(u_t)_+:=\max\{u_t,0\}$.
The fact that the function $f(x):=-xe^{-x}$ is bounded over the interval $[0,+\infty)$ implies the uniform estimate $a_t\geq -C$. We prove now the upper bound of $a_t$.  
Perelman show this by using the monotonicity of his $\mu$ functional along the K\"ahler-Ricci flow. We realize that the the upper bound of $a_t$ follows in a classical way by using the generalized Bochner-Kodaira Formula. In fact using the K\"ahler-Ricci flow identity \eqref{EvolUfunc} and the identity $\Delta_te^{-u_t}=(2|\partial u_t|^2_t-\Delta_tu_t)e^{-u_t}$, we find
\begin{eqnarray*}
-\dot{a}_t&=&\dashint\limits_X \Big[\dot{u}_t-u_t\left(\dot{u}_t-\Delta_tu_t/2\right)\Big]e^{-u_t}\omega^n_t
=
\dashint\limits_X\Big[\dot{u}_t-u_t\left(u_t+a_t\right)\Big]e^{-u_t}\omega^n_t
\\
&=& 
\frac{1}{2}\,\dashint\limits_X\Delta_tu_t e^{-u_t}\omega^n_t+\dashint\limits_X u_te^{-u_t}\omega^n_t+a_t-\dashint\limits_X u^2_te^{-u_t}\omega^n_t+a_t^2
\\
&=&
\dashint\limits_X |\partial u_t|^2_te^{-u_t}\omega^n_t-\dashint\limits_X u^2_te^{-u_t}\omega^n_t+a_t^2
\end{eqnarray*}
By the Poincarr\'e type inequality \ref{PoincIneq} in the Fano case we deduce $-\dot{a}_t\geq 0$. This implies the upper bound of the normalizing constants $a_t$. 
%(remark that $\dot{a}_{t_0}=0$ for some time $t_0\geq 0$, if and only if $\dot{a}_t=0$ for all $t\geq 0$, if and only if $\omega_t$ is a Gradient-K\"ahler-Ricci soliton.)
Consider now Perelman's functional
\begin{eqnarray*}
{\cal W}(\omega,f,\tau):=(4\pi c_1\tau)^{-n}\int\limits_X\Big[\tau\left(|\nabla_{\omega}\, f|^2_{\omega}+\Scal_{\omega}\right)+f-2n\Big]e^{-f}\omega^n.
\end{eqnarray*}
Using the identity $\Delta_tu_t=2n-\Scal_t$ we get
\begin{eqnarray*}
{\cal W}_t:={\cal W}(\omega_t,u_t,1/2)
&=&
\dashint\limits_X\left[\frac{1}{2}\left(|\nabla_t\,u_t|^2_t+\Scal_t\right)+u_t-2n\right]e^{-u_t}\omega^n_t
\\
&=& 
\dashint\limits_X\left[\frac{1}{2}\left(|\nabla_t\,u_t|^2_t-\Delta_tu_t\right)+u_t\right]e^{-u_t}\omega^n_t-n
\\
&=& 
\dashint\limits_X\left(\frac{1}{2}\Delta_te^{-u_t}+u_te^{-u_t}\right)\omega^n_t-n
\\
&=& -a_t-n\,.
%\dashint\limits_X u_te^{-u_t}\omega^n_t=-a_t,
\end{eqnarray*}
So we find the inequality $\dot{{\cal W}}_t=-\dot{a}_t\geq 0$. Suppose now that $\dot{{\cal W}}_t=0$ for some time $t$. Then we get the equality case in the Poincarr\'e inequality
$$
\int\limits_X|\nabla_t\,\theta_t|^2_t\,e^{-u_t}\omega^n_t=2\int\limits_X\theta_t^2\,e^{-u_t}\omega^n_t\,,
$$
with $\theta_t:=u_t-\dashint_X u_te^{-u_t}\omega^n_t=u_t+a_t$.
The variational characterization of the first non zero eigenvalue $\lambda_1(\hat\Delta_t)$ of the generalized Laplacian  $\hat\Delta_t:=\Delta_{t,-u_t}$ implies $2\geq\lambda_1(\hat\Delta_t)$. By the other hand the generalized Bochner-Kodaira formula implies $\lambda_1(\hat\Delta_t)\geq 2$. So $\lambda_1(\hat\Delta_t)=2$ and $\hat\Delta_tu_t+2(u_t+a_t)=0$. By plugging this in to the generalized Bochner-Kodaira formula we get
$$
\int\limits_X|\bar{\partial}\,\nabla^{1,0}_t u_t |^2_t\,e^{-u_t}\omega^n_t
=
-\int\limits_X\left<\partial\hat\Delta_tu_t ,\partial u\right>_te^{-u_t}\omega^n_t-\int\limits_X|\nabla_t\,u_t|^2_t\,e^{-u_t}\omega^n_t=0\,.
$$
So $\omega_t$ is a K\"ahler-Ricci soliton and this will hold for all times. We have prove in conclusion the  proposition \ref{MonWkrfRP}.
\\
\\
{\bf The generalized functionals by Aubin.}
\\
The generalized functionals
$I_{\omega},\,J_{\omega}:{\cal P}_{\omega}\rightarrow [0,+\infty)$ by Aubin, \cite{Aub1} are defined by the formulas
\begin{eqnarray*}
I_{\omega}(\varphi )&:=&\dashint\limits_X\varphi \left(\omega ^n-\omega^n_{\varphi}\right)
=\sum_{k=0}^{n-1} \,
\dashint\limits_X i\partial \varphi \wedge \bar{\partial}\varphi \wedge\omega ^k\wedge \omega^{n-k-1}_{\varphi}
\\
J_{\omega}(\varphi )&:=&\sum_{k=0}^{n-1} \frac{k+1}{n+1}\,
\dashint\limits_X i\partial \varphi \wedge \bar{\partial}\varphi \wedge\omega ^k\wedge \omega^{n-k-1}_{\varphi}
\\
&=&\dashint\limits_X\varphi \,\omega ^n-\frac{1}{n+1}\sum_{k=0}^n\,\dashint\limits_X\varphi \,\omega ^k\wedge \omega ^{n-k}_{\varphi}.
\end{eqnarray*}
We have the obvious inequalities, $0\leq I_{\omega}\leq (n+1)J_{\omega}$. 
\\
\\
{\bf The K-energy functional of the anticanonical class $2\pi c_1$.}
We remind that the Einstein equation 
$
\Ric(\omega _{\varphi} )=\omega _{\varphi}, 
$
is equivalent to the constant scalar curvature equation $\Scal(\omega _{\varphi})=2n$.
This last equation is the  Euler-Lagrange equation of Mabuchi's \cite{Mab} K-energy functional $\nu_{\omega}:{\cal P}_{\omega}\rightarrow \R$
$$
\nu _{\omega}(\varphi):=
\dashint\limits_X\left(\log\,\frac{\omega _{\varphi} ^n}{\omega ^n}+\varphi-h_{\omega}  \right)\omega^n_{\varphi}
-
\frac{1}{n+1}\sum_{k=0}^n\,\dashint\limits_X\varphi \,\omega ^k\wedge \omega ^{n-k}_{\varphi}
+ \dashint\limits_X h_{\omega}\omega ^n.
$$
In fact for every $\ci$ path $(\varphi _t)_{t\in (-\varepsilon ,\varepsilon )} \subset{\cal P}_{\omega}$ we have the identity
\begin{eqnarray}\label{der-Kenerg} 
\frac{d}{dt}\nu_{\omega}(\varphi_t)=-\frac{1}{2}  \,\dashint\limits_X
\dot{\varphi}_t\Big( \Scal(\omega _t )-2n\Big)\omega ^n_t\,, 
\end{eqnarray}  
where $\dot{\varphi}_t:=\frac{\partial}{\partial t}\varphi _t$ and $\omega_t:=\omega_{\varphi_t}$. 
We remark that under the  K\"ahler-Ricci flow we have the identity
$
\Scal(\omega _t)=2n-\Delta _{\omega_t}\dot{\varphi}_t
$
.
Then using the identity \eqref{der-Kenerg} we deduce the inequality
$$
\frac{d}{dt}\nu _{\omega}(\varphi _t)=2^{-1} \dashint\limits_X\dot{\varphi}_t\Delta _{\omega_t}\dot{\varphi}_t\,\omega ^n_t=-n\dashint\limits_X
i\partial \dot{\varphi}_t\wedge \bar{\partial}\dot{\varphi}_t\wedge \omega ^{n-1}_t \leq 0,
$$
which shows that the K-energy decreases under the K\"ahler-Ricci flow. We remind also the following Tian's \cite{Tia} fundamental result.
\begin{theorem}{\bf (Tian's $G$-properness)}
Let $X$ be a Fano manifold admitting a $G$-invariant K\"ahler-Einstein metric $\hat\omega\in 2\pi c_1$. Then there exists two constants $\delta>0,\,C>0$ such that the inequality $\nu_{\hat\omega}(\varphi)\geq J_{\hat\omega}(\varphi)^{\delta}-C$ hold for all $G$-invariant potentials $\varphi\in {\cal P}_{\hat\omega}$.
\end{theorem}
By using the cocycle condition we deduce that for all $G$-invariant K\"ahler metrics $\omega\in 2\pi c_1$ there exist an increasing function $\mu:\R\rightarrow [c,+\infty)$, with 
$\lim_{t\rightarrow +\infty}\mu(t)$ 
$=+\infty$ such that $\nu_{\omega}(\varphi)\geq \mu( J_{\omega}(\varphi))$ for all $G$-invariant potentials $\varphi\in {\cal P}_{\omega}$.
%%%%%%%%%%%%%%%%%%%%%%%%%%%%%%%%%%%%%%%%%%%%%%%%%%%%%%%%%%%%%%%%%%%%%%%%%%%%%%%%%%%%%%%%%%%%%%%%%%%%%%%%%%%%%%%%%%%%%%%%%%%%%%%%%%%%%%%%%%%%%%%%%%%%%%%%%%%%%%%%%%%%%%%%%%%%%%%%%%%%%%%%%%%%%%%%%%%%%%%%%%%%%%%%%%%%%%%%%%%%%%%%%%%%%%%%%%%%%%%%%%%%%%%%%%%%%%%%%%%%%
%%%%%%%%%%%%%%%%%%%%%%%%%%%%%%%%%%%%%%%%%%%%%%%%%%%%%%%%%%%%%%%%%%%%%%%%%%%%%%%%%%%%%%%%%%%%%%%%%%%%%%%%%%%%%%%%%%%%%%%%%%%%%%%%%%%%%%%%%%%%%%%%%%%%%%%%%%%%%%%%%%%%%%%%%%%%%%%%%%%%%%%%%%%%%%%%%%%%%%%%%%%%%%%%%%%%%%%%%%%%%%%%%%%%%%%%%%%%%%%%%%%%%%%%%%%%%%%%%%%%%%%%%%%%%%%%
\section{Tian-Zhu's $C^0$-uniform estimate}
%%%%%%%%%%%%%%%%%%%%%%%%%%%%%%%%%%%%%%%%%%%%%%%%%%%%%%%%%%%%%%%%%%%%%%%%%%%%%%%%%%%%%%%%%%%%%%%%%%%%%%%%%%%%
%%%%%%%%%%%%%%%%%%%%%%%%%%%%%%%%%%%%%%%%%%%%%%%%%%%%%%%%%%%%%%%%%%%%%%%%%%%%%%%%%%%%%%%%%%%%%%%%%%%%%%%%%%%%
We start by proving the following elementary lemma
\begin{lem}\label{DsVol-J}
Let $(X,\omega)$ be a polarized Fano manifold with $\omega\in 2\pi c_1$, let $(\omega_t)_t$ be a K\"ahler-Ricci flow and let and $G$ be a compact maximal subgroup of the identity component of the group of automorphisms of $X$.
\\
{\bf A)}. Suppose there exist a constant $k>0$ such that $\omega^n_t\geq k\omega^n$ for all times $t\geq 0$. Then the Aubin's Functional $J_{\omega}$ is uniformly bounded along this K\"ahler-Ricci flow.
\\
{\bf B)}. Suppose $X$ admits a $G$-invariant K\"ahler-Einstein metric. Then the K\"ahler-Ricci flow with $G$-invariant initial metric $\omega$ satisfies the uniform estimate $\omega_t^n\geq k\,\omega^n$, $k>0$ for all times $t\geq 0$.
\end{lem}
So part B of this lemma prove one implication in theorem 1.
\\
\\
$Proof$. Set $\hat\varphi_t:=\varphi_t+c_t$. By writing the K\"ahler-Ricci flow equation under the form
$$
\hat\varphi_t=\dot{\varphi}_t-\log\,\frac{\omega _t ^n}{\omega ^n}+h_{\omega}\,,
$$
using the Perelman's uniform estimate $|\dot{\varphi}_t|\leq C$ and the inequality $\omega_t^n\geq k\omega^n$ we find the uniform estimate
$\hat\varphi_t\leq C-\log k+h_{\omega}$.
Reminding the expression of the K-energy functional we deduce the identity along the K\"ahler-Ricci flow
\begin{eqnarray}\label{KenKRfloNorm}
\nu_{\omega}(\varphi_t)=\dashint\limits_X\dot{\varphi}_t\,\omega^n_t+J_{\omega}(\varphi_t)-\dashint\limits_X \hat\varphi_t\,\omega^n
+\dashint\limits_X h_{\omega}\,\omega^n.
\end{eqnarray}
Then using; the fact that the K-energy functional is nonincreasing along the K\"ahler-Ricci flow, the inequality $-\dashint_X\dot{\varphi}_t\,\omega^n_t\leq 0$ (which follows from the integral normalization of $\dot{\varphi}_t$) and the previous estimate $\hat\varphi_t\leq C$, we deduce the uniform estimate 
$$
0\leq J_{\omega}(\varphi_t)=\nu_{\omega}(\varphi_t)-\dashint\limits_X\dot{\varphi}_t\,\omega^n_t+\dashint\limits_X \hat\varphi_t\,\omega^n
-\dashint\limits_X h_{\omega}\,\omega^n\leq \nu_{\omega}(\varphi_0)+C.
$$
%Using the fact that $J_{\omega}\geq 0$ in the identity \eqref{KenKRfloNorm} we have
%$$\nu_{\omega}(\varphi_0)\geq \nu_{\omega}(\varphi_t)=\dashint\limits_X\dot{\varphi}_t\,\omega^n_t-\dashint\limits_X(\varphi_t+c_t)\,\omega^n+\dashint\limits_X h_{\omega}\,\omega^n,$$
%which implies the uniform estimate
%$$-\dashint\limits_X(\varphi_t+c_t)\,\omega^n\leq C.$$
We prove now part B. The existence of a $G$-invariant K\"ahler-Einstein metric implies the $G$-properness of the K-energy functional. Then using the fact that the K-energy functional is nonincreasing along the K\"ahler-Ricci flow we deduce that the energy functional $J_{\omega}$ is bounded along the K\"ahler-Ricci flow $(\omega_t)_t$ with $G$-invariant initial metric $\omega$. Then the identity \eqref{KenKRfloNorm} combined with the fact that the K-energy functional is bounded from below implies the uniform estimate $\dashint_X\hat\varphi_t\,\omega^n\leq C$.
By the properties of the Green function we deduce the inequality
\begin{eqnarray*}
\hat\varphi_t\leq \dashint\limits_X\hat\varphi_t\,\omega^n+C'\leq C\,.
\end{eqnarray*}
This uniform estimate is equivalent to the uniform estimate $\omega_t^n\geq k\,\omega^n$ by means of the K\"ahler-Ricci flow equation and Perelman's uniform estimate $|\dot{\varphi}_t|\leq C$. \hfill $\Box$
\\
We remind now the following result \cite{Ti-Zh}.
\begin{prop}\label{Kolo}
Let $(X,\omega)$ be a compact K\"ahler manifold of complex dimension $n$, let $\varphi\in {\cal P}_{\omega}$ and $f:=\omega_{\varphi}^n/\omega^n$. Then for all $\varepsilon\in (0,\varepsilon_0]$, $\delta\in (0,\delta_0]$ there exists constants $C,C'>0$ depending only on $\omega, \,\varepsilon_0,\,\delta_0$ such that 
$$
\Osc(\varphi)\leq C\left(\frac{1}{\varepsilon\delta}\right)^{n+\delta}\|f\|^{\delta}_{L^{1+\varepsilon}(X,\omega)}+C'\,.
$$
\end{prop}
\begin{prop}
Let $(X,\omega)$ be a polarized Fano manifold with $\omega\in 2\pi c_1$ and let $(\omega_t)_t$ be a K\"ahler-Ricci flow  admiting a constant $k>0$ such that $\omega^n_t\geq k\omega^n$ for all times $t\geq 0$. Then this K\"ahler-Ricci flow satisfies the uniform estimate
$|\varphi_t+c_t|\leq K_0$, for some constant $K_0>0$ independent of $t\geq 0$.
\end{prop}
$Proof$. The argument here is the same as in \cite{Ti-Zh}.  We start proving the uniform estimate $|\max_X\hat\varphi_t|\leq C$ for all $t\geq 0$. We first remark that if a real function $u$ satisfies the integral equality $\int_X(e^{-u}-1)\omega^n=0$ then $\max_Xu\geq0$. 
In fact if not $0>\max_Xu\geq u$ and this implies $e^{-u}>1$, which contradict the integral equality. By definition of K\"ahler-Ricci flow we have the integral identity 
\begin{eqnarray}\label{IntKRFC0}
\int\limits_Xe^{h-\hat\varphi_t}\omega^n=\int\limits_Xe^{-\dot{\varphi}_t}\omega^n_t=\int\limits_X\omega^n,
\end{eqnarray}
for all $t\in [0,+\infty)$. Then applying the previous remark with $u:=\hat\varphi_t-h$, we find the inequality
$\max_X(\hat\varphi_t-h)\geq 0$, which gives
the estimate $\max_X\hat\varphi_t\geq -C$. Moreover the argument in the proof of A of lemma \ref{DsVol-J} implies $\hat\varphi_t\leq C$. The equality \eqref{IntKRFC0} implies that the function $h-\hat\varphi_t$ change signs and so we get $\|h-\hat\varphi_t\|_{C^0(X)}\leq \Osc(h-\hat\varphi_t)$, which implies $\|\hat\varphi_t\|_{C^0(X)}\leq \Osc(\hat\varphi_t)+C$. By proposition \ref{Kolo} we need to prove a uniform bound for the integral $\int_Xe^{-(1+\varepsilon)\hat\varphi_t}\omega^n$ for some $\varepsilon>0$. Set $\theta_t:=\max_X\varphi_t-\varphi_t\geq 0$. Then
\begin{eqnarray*}
e^{-(1+\varepsilon)\hat\varphi_t}\omega^n=e^{\varepsilon\theta_t-\varepsilon\max\hat\varphi_t-\hat\varphi_t}
\omega^n\leq Ce^{\varepsilon\theta_t-\hat\varphi_t}
\omega^n\leq C'e^{\varepsilon\theta_t}
\omega^n_t\,.
\end{eqnarray*}
The last inequality follows from the K\"ahler-Ricci flow equation and Perelman's uniform estimate $|\dot{\varphi}_t|\leq C$. So it is sufficient to prove an uniform bound for the integral $\int_Xe^{\varepsilon\theta_t}
\omega^n_t$. In order to prove this we consider the classic inequality
$$
0\leq I_{\omega}(\varphi_t)=
\dashint\limits_X\hat\varphi_t
\left(\omega^n-\omega^n_t\right)\leq (n+1)J_{\omega}(\varphi_t)\leq C\,,
$$ 
where $C>0$ is the uniform constant provided by lemma \ref{DsVol-J}. We deduce
$$
-\dashint\limits_X\hat\varphi_t\,\omega^n_t\leq C-\dashint\limits_X\hat\varphi_t\,\omega^n\leq 2C.
$$
The last inequality follows from the estimate $\dashint_Xe^{-\hat\varphi_t}\,\omega^n\leq C$ that we get from the identity \eqref{IntKRFC0}. So we have obtain the uniform estimate
\begin{eqnarray}\label{tetaC0}
0\leq \dashint\limits_X\theta_t\,\omega^n_t\leq C
\end{eqnarray}
%As in \cite{Ti-Zh} we consider 
For all integers $p\geq 1$ we have the equalities
\begin{eqnarray*}
\int\limits_X\theta^p_t\left(\omega^n_t-\omega^{n-1}_t\wedge \omega\right)
&=&-\int\limits_X\theta^p_t \,i\partial\bar{\partial}\theta_t\wedge \omega^{n-1}_t
\\
&=&
p\int\limits_X\theta^{p-1}_t i\partial \theta_t\wedge \bar{\partial}\theta_t\wedge \omega^{n-1}_t
\\
&=&p\int\limits_X\theta^{\frac{p-1}{2}}_t i\partial \theta_t\wedge \theta^{\frac{p-1}{2}}_t\bar{\partial}\theta_t\wedge \omega^{n-1}_t
\\
&=&\frac{4p}{(p+1)^2}\int\limits_Xi\partial \theta^{\frac{p+1}{2}}_t  \wedge \bar{\partial}\theta^{\frac{p+1}{2}}_t\wedge \omega^{n-1}_t
\\
&=&\frac{4p}{n(p+1)^2}\int\limits_X |\partial \theta^{\frac{p+1}{2}}_t|^2_t\,\omega^n_t\,.
\end{eqnarray*}
This implies the inequality
\begin{eqnarray}\label{IdiotC0Inte}
\int\limits_X |\partial \theta^{\frac{p+1}{2}}_t|^2_t\,\omega^n_t\leq \frac{n(p+1)^2}{4p}\int\limits_X\theta^p_t\,\omega^n_t.
\end{eqnarray}
Remember now the K\"ahler-Ricci flow identity $i\partial\bar{\partial}\dot{\varphi}_t=\omega_t-\Ric(\omega_t)$.
Then as in \cite{Ti-Zh} by applying the Poincarr\'e type inequality (corollary \ref{PoincIneq}) to the function $\theta^{\frac{p+1}{2}}_t$, with metric $\omega_t$ and $h=-\dot{\varphi}_t$, we deduce
$$
\int\limits_X |\partial \theta^{\frac{p+1}{2}}_t|^2_t\,e^{-\dot{\varphi}_t}\omega^n_t\geq 
\int\limits_X\theta^{p+1}_te^{-\dot{\varphi}_t}\omega^n_t-(2\pi c_1)^{-n}\left(\;\int\limits_X\theta^{\frac{p+1}{2}}_te^{-\dot{\varphi}_t}\omega^n_t\right)^2.
$$
By applying the H\"older inequality to the last therm, using Perelman uniform estimate $|\dot{\varphi}_t|\leq C$ and the inequality
\eqref{IdiotC0Inte} we deduce
$$
\int\limits_X\theta^{p+1}_te^{-\dot{\varphi}_t}\omega^n_t\leq Cp\int\limits_X\theta^p_t\,\omega^n_t+C\int\limits_X\theta^p_t\,e^{-\dot{\varphi}_t}\omega^n_t\cdot \int\limits_X\theta_t\,e^{-\dot{\varphi}_t}\omega^n_t.
$$
Using again the estimate $|\dot{\varphi}_t|\leq C$ and the estimate \eqref{tetaC0} we find
$$
\int\limits_X\theta^{p+1}_t\,\omega^n_t\leq C(p+1)\int\limits_X\theta^p_t\,\omega^n_t.
$$
By iteration $\int_X\theta^p_t\,\omega^n_t\leq C^pp!$. Thus
$$
\int\limits_Xe^{\varepsilon\theta_t}\omega^n_t=\sum_{p=0}^{\infty}\frac{\varepsilon^p}{p!}\int\limits_X\theta^p_t\,\omega^n_t\leq \sum_{p=0}^{\infty}(\varepsilon C)^p\,.
$$
So we choose $0<\varepsilon<1/C$.\hfill$\Box$
\newpage
%%%%%%%%%%%%%%%%%%%%%%%%%%%%%%%%%%%%%%%%%%%%%%%%%%%%%%%%%%%%%%%%%%%%%%%%%%%%%%%%%%%%%%%%%%%%%%%%%%%%%%%%%%%%%%%%%%%%%%%%%%%%%%%%%%%%%%%%%%%%%%%%%%%%%%%%%%%%%%%%%%%%%%%%%%%%%%%%%%%%%%%%%%%%%%%%%%%%%%%%%%%%%%%%%%%%%%%%%%%%%%%%%%%%%%%%%%%%%%%%%%%%%%%%%%%%%%%%%%%%%%%%%%%%
\section{The Yau's $C^2$ and Calabi's $C^3$ uniform estimates for the K\"ahler-Ricci flow}
%%%%%%%%%%%%%%%%%%%%%%%%%%%%%%%%%%%%%%%%%%%%%%%%%%%%%%%%%%%%%%%%%%%%%%%%%%%%%%%%%%%%%%%%%%%%%%%%%%%%%%%%%%%%%%%%%%%%%%%%%%%%%%%%%%%%%%%%%%%%%%%%%%%%%%%%%%%%%%%%%%%%%%%%%%%%%%%%%%%%%%%%%%%%%%%%%%%%%%%%%%%%%%%%%%%%%%%%%%%%%%%%%%%%%%%%%%%%%%%%%%%%%%%%%%%%%%%%%%%%%%%%%%%%%%
%\subsection{Yau's $C^2$-uniform estimate}
%We will give first a version of the proof of the Yau's $C^2$-uniform estimate for the complex Monge-Amp\`ere equation \cite{Yau}. 
We start with some notations and definitions. Let $(X,\omega )$ be a K\"ahler manifold of complex dimension $n$. 
Consider the function $\lambda _1^{\omega}:X\rightarrow \R$ defined by the formula
$$
\lambda _1^{\omega}(x):=\min_{\xi \in T^{\otimes 2}_{X,x}\smallsetminus 0_x } 
{\cal C}^{\omega}_{_{X,J}}(\xi ,\xi)|\xi|_{\omega} ^{-2}. 
$$
So $\lambda _1^{\omega}(x)$ is the smallest eigenvalue of the Chern Curvature form ${\cal C}^{\omega}_{_{X,J}}(x)$. 
It is well known (see \cite{Kat}, chap II, sec 5.1, theorem 5.1, pag 107) that the function  $\lambda _1^{\omega}$ is continuous.
In order to simplify the notations we will use Einstein convention on sums. Moreover we will note by $\Tr_{\varphi}$ the trace operator corresponding to the metric $\omega _{\varphi}$. With such notations we have the following proposition which is obtained by some computations in \cite{Yau}.
\begin{prop}\label{RiccCaYau} 
Let $(X,\omega )$ be a K\"ahler manifold of complex dimension $n$. Then for every potential $\varphi \in {\cal P}_{\omega}$ for the K\"ahler metric $\omega$ we have the intrinsic inequality
$$
2\Tr_{\omega} \Ric(\omega_{\varphi})\geq  
-\Delta _{\varphi}\Delta _{\omega}\varphi 
+4\lambda _1^{\omega}\,(2n +\Delta_{\omega }\varphi )\Tr_{\varphi}\omega
+\frac{2|\partial \Delta _{\omega}\varphi  |^2_{\varphi }}{2n +\Delta_{\omega }\varphi}.
$$  
\end{prop}
$Proof$. Let $(z_1,...,z_n)$ be $\omega$-geodesic holomorphic coordinates with center a point $x$ such that the metric $\omega _{\varphi}$ can be written in diagonal form in $x$.  Explicitly 
$\omega =\frac{i}{2}\omega _{l,\bar{r}}\,dz_l\wedge d\bar{z}_r$ and 
$\omega_{\varphi}  =\frac{i}{2}(\omega_{\varphi})_{l,\bar{r}}\,dz_l\wedge d\bar{z}_r$, with
\begin{eqnarray*}
&\displaystyle{
\omega _{l,\bar{r}}=\delta _{l,r}-C^{j,\bar{k}}_{r,l}z_j\bar{z}_k+O(|z|^3), }&
\\
&\displaystyle{
{\cal C}_{\omega }(T_{_{X,J} } )(x) =
 C^{j,\bar{k}}_{l,r}\,(dz_j\wedge d\bar{z}_k)\otimes dz_r\otimes_{_{J}}\frac{\partial }{\partial z_l},}&
\\
&\displaystyle{
{\cal C}^{\omega}_{_{X,J}}(x)=C^{j,\bar{k}}_{r,l}\,dz_j\otimes dz_l\otimes d\bar{z}_k\otimes d\bar{z}_r,}&
\\
&\displaystyle{
\overline{C^{j,\bar{k}}_{l,r}}=C^{k,\bar{j}}_{r,l},\;C^{j,\bar{k}}_{l,r}=C^{r,\bar{k}}_{l,j}=C^{j,\bar{l}}_{k,r},}&
\\
&\displaystyle{
(\omega_{\varphi})_{l,\bar{r}} =\delta _{l,r}+2\varphi_{l,\bar{r}}+O(|z|),\quad 2\varphi_{l,\bar{r}}(0)= 2\delta _{l,r} \varphi_{l,\bar{l}}(0)>-1,}&
\end{eqnarray*}
where $\varphi_{l,\bar{r}}:=\frac{\partial^2\varphi}{\partial z_l\partial\bar{z}_r }$.
In particular we deduce the following expressions for the inverse matrixs
\begin{eqnarray*}
%&\displaystyle{
\omega^{l,\bar{r}}=\delta _{l,r}+C^{j,\bar{k}}_{r,l}z_j\bar{z}_k+O(|z|^3), %}&
\qquad
%&\displaystyle{
(\omega_{\varphi})^{l,\bar{r}}=\frac{\delta _{l,r}}{1+2\varphi_{l,\bar{l}} }+O(|z|).%}&  
\end{eqnarray*}
Moreover we deduce the local expressions
\begin{eqnarray*}
&\displaystyle{
\Delta _{\omega}\varphi =4\omega^{l,\bar{r}}\varphi_{r,\bar{l}}=4(\varphi_{l,\bar{l}}
+C^{j,\bar{k}}_{r,l}\varphi_{r,\bar{l}}\,z_j\bar{z}_k  ) +O(|z|^3),}&
\\
&\displaystyle{
\Delta _{\varphi}u=4(\omega_{\varphi})^{l,\bar{r}}u_{r,\bar{l}}=
\frac{4u_{l,\bar{l}}}{1+2\varphi_{l,\bar{l}} }+O(|z|)}&
\end{eqnarray*}
for every smooth function $u$. Using this two expressions we find the equality at the point $x$
\begin{eqnarray}\label{DelDel} 
\Delta _{\varphi}\Delta _{\omega}\varphi
=
\frac{4^2}{1+2\varphi_{l,\bar{l}} }
(\varphi_{l,\bar{l},k,\bar{k}}+C^{l,\bar{l}}_{r,k}\varphi_{r,\bar{k}})
=
\frac{4^2}{1+2\varphi_{l,\bar{l}} }
(\varphi_{l,\bar{l},k,\bar{k}}+C^{l,\bar{l}}_{k,k}\varphi_{k,\bar{k}})\,.
\end{eqnarray}
Using the expression \eqref{CoefCurv} of the coefficients  of the curvature form respect to the complex frame $(\zeta_k):=(\partial/\partial z_k)$, we find the following local expression for the Ricci tensor
\begin{eqnarray*}
\Ric(\omega)=-\Big(i\partial\bar{\partial}\omega_{l,\bar{r}}-i\partial\omega_{l,\bar{s} }\,\omega ^{s,\bar{t}}\wedge\bar{\partial}\omega_{t,\bar{r}}\Big)\omega ^{r,\bar{l}}.  
\end{eqnarray*}
All the computations that will follow are refereed to the point $x$.
Expanding the analogue expression for $\Ric(\omega_{\varphi})$ we get the equality.
\begin{eqnarray*}
\Ric(\omega_{\varphi})= 
\Big(iC_{r,l}-2i\partial\bar{\partial}\varphi_{l,\bar{r}} +4i\partial\varphi_{l,\bar{s} }\,\omega_{\varphi} ^{s,\bar{t}}\wedge\bar{\partial}\varphi_{t,\bar{r}}\Big)
\omega_{\varphi} ^{r,\bar{l}}.  
\end{eqnarray*}
Taking the trace respect to $\omega$ of the Ricci tensor $\Ric(\omega_{\varphi})$ we find the expression
\begin{eqnarray*}
\Tr_{\omega} \Ric(\omega_{\varphi})&=&
4\Big(C^{k,\bar{k}}_{r,l}-2\varphi_{k,\bar{k}, l,\bar{r}} +4\varphi_{k,l,\bar{s} }\,\omega_{\varphi} ^{s,\bar{t}}\,\varphi_{\bar{k}, t,\bar{r}}\Big)\omega_{\varphi} ^{r,\bar{l}}    
\\
&=& 4\Big(C^{k,\bar{k}}_{l,l}-
2\varphi_{k,\bar{k}, l,\bar{l}}+
\frac{4\varphi_{k,l,\bar{s} }\,\varphi_{\bar{k},s,\bar{l} }}{1+2\varphi_{s,\bar{s} }} \Big)\frac{1}{1+2\varphi_{l,\bar{l}}}.    
\end{eqnarray*}
Using the symmetry $C^{k,\bar{k}}_{l,l}=C^{l,\bar{l}}_{k,k}\in \R$ and the identity \eqref{DelDel} we find the equality
\begin{eqnarray}\label{TrRic}
\Tr_{\omega} \Ric(\omega_{\varphi})=-\frac{1}{2}\Delta _{\varphi}\Delta _{\omega}\varphi 
+
4C^{l,\bar{l}}_{k,k}\frac{1+2\varphi _{k,\bar{k}}}{1+2\varphi_{l,\bar{l}}}+
\frac{16\varphi_{k,l,\bar{s} }\,\varphi_{\bar{k},s,\bar{l} }}{(1+2\varphi_{s,\bar{s}})
(1+2\varphi_{l,\bar{l}})}.
\end{eqnarray}
The inequality
$
C^{l,\bar{l}}_{k,k}={\cal C}^{\omega}_{_{X,J}}
(\frac{\partial}{\partial x_l}\otimes \frac{\partial }{\partial x_k},\frac{\partial}{\partial x_l}\otimes \frac{\partial }{\partial x_k})(x)\geq \lambda _1^{\omega}(x),  
$
implies the inequality
$$
4C^{l,\bar{l}}_{k,k}\frac{1+2\varphi _{k,\bar{k}}}{1+2\varphi_{l,\bar{l}}}\geq 
2\lambda _1^{\omega}\,(2n +\Delta_{\omega }\varphi )\Tr_{\varphi}\omega\,(x).%%%%%%%%%
$$
Then the conclusion of the proof of the proposition will follows from the inequality
\begin{eqnarray}\label{ineqTirdDer} 
\frac{16\varphi_{k,l,\bar{s} }\,\varphi_{\bar{k},s,\bar{l} }}{(1+2\varphi_{s,\bar{s} })
(1+2\varphi_{l,\bar{l}})}
\geq
\frac{|\partial \Delta _{\omega}\varphi  |^2_{\varphi }}{2n +\Delta_{\omega }\varphi}
.   
\end{eqnarray}
Let prove this inequality. We have
\begin{eqnarray*}
|\partial \Delta _{\omega}\varphi  |^2_{\varphi }
&=&
2\sum_{k,l} \omega_{\varphi}^{k,\bar{l}}\,\partial_l\Delta_{\omega}  \varphi 
\,\partial_{\bar{k}}\Delta _{\omega}\varphi
=
\sum_{j,k,l}
\frac{2\cdot 4^2\varphi_{l,j,\bar{j} }\,\varphi_{\bar{l},k,\bar{k} }}{1+2\varphi_{l,\bar{l}}}
\\
&=&
\sum_l\frac{2\cdot 4^2}{1+2\varphi_{l,\bar{l}}}  \left|
\sum_j
\frac{\varphi_{l,j,\bar{j} }}{\sqrt{1+2\varphi_{j,\bar{j}}} }
\sqrt{1+2\varphi_{j,\bar{j}}}\, \right|^2.
\end{eqnarray*}
Applying the Cauchy-Schwartz inequality to the norm, we find the inequality
\begin{eqnarray*}
|\partial \Delta _{\omega}\varphi  |^2_{\varphi }
&\leq&
\sum_l\frac{2\cdot 4^2}{1+2\varphi_{l,\bar{l}}} 
\left(
\sum_j
\frac{|\varphi_{l,j,\bar{j} }|^2}{1+2\varphi_{j,\bar{j}} }\right) 
\left(
\sum_k
(1+2\varphi_{k,\bar{k}})\right)
\\
&=&
(2n+\Delta _{\omega}\varphi)\sum_{j,l} 
\frac{16|\varphi_{l,j,\bar{j} }|^2}{(1+2\varphi_{j,\bar{j}})(1+2\varphi_{l,\bar{l}}) }
\\
&\leq&
(2n+\Delta _{\omega}\varphi)\sum_{k,l,j}\frac{16\varphi_{k,l,\bar{j} }\,\varphi_{\bar{k},j,\bar{l} }}{(1+2\varphi_{j,\bar{j} })
(1+2\varphi_{l,\bar{l}})}\,,
\end{eqnarray*}
which conclude the proof of the inequality \eqref{ineqTirdDer}.\hfill$\Box$
\\
\\
%%%%%%%%%%%%%%%%%%%%%%%%%%%%%%%%%%%%%%%%%%%%%%%%%%%%%%%%%%%%%%%%%%%%%%%%%%%%%%%%%%%%%%%%%%%%%%%%%%%%%%%%%%%%%%%%%%%%%%%%%%%%%%%%%%%%%%%%%%%%%%%%%%%%%%%%%%%%%%%%%%%%%%%%%%%%%%%%%%%%%%%%%%%%%%%%%%%%%%%%%%%%%
%                                    CALABI   C3-ESTIMATE
%%%%%%%%%%%%%%%%%%%%%%%%%%%%%%%%%%%%%%%%%%%%%%%%%%%%%%%%%%%%%%%%%%%%%%%%%%%%%%%%%%%%%%%%%%%%%%%%%%%%%%%%%%%%%%%%%%%%%%%%%%%%%%%%%%%%%%%%%%%%%%%%%%%%%%%%%%%%%%%%%%%%%%%%%%%%%%%%%%%%%%%%%%%%%%%%%%%%%%%%%%%%
%\subsection{Calabi's computation for the $C^3$ uniform estimate}
%Let $(X,\omega )$ be a K\"ahler manifold and let ${\varphi\in\cal P}_{\omega}$ be a potential of $\omega$.
We will note by $\omega^*$ and $\omega^*_{\varphi}$ the corresponding dual elements of $\omega$ and $\omega_{\varphi}$. Let $h^*$ and $h^*_{\varphi}$ the corresponding hermitian metrics over the complex vector bundle $T^*_{_{X,J}}$. In local complex coordinates we have the expressions $\omega^*=2i\,\omega^{l\bar{k}}\frac{\partial }{\partial z_k}\wedge\frac{\partial }{\partial \bar{z}_l}$ and $h^*=4\,\omega^{l\bar{k}}\frac{\partial }{\partial z_k}\otimes \frac{\partial }{\partial \bar{z}_l}$ .
%$$\omega^*=2i\,\omega^{l\bar{k}}\frac{\partial }{\partial z_k}\wedge\frac{\partial }{\partial \bar{z}_l}
%\quad\mbox{and}\quad h^*=4\,\omega^{l\bar{k}}\frac{\partial }{\partial z_k}\otimes \frac{\partial }{\partial \bar{z}_l}.$$
We remind also that if $(V,J)$ is a complex vector space equipped with a hermitian metric $h$ then the corresponding hermitian metric $h_{_\C}$ over the complexified vector space $(V\otimes_{_\R}\C, i)$ is defined by the formula
$$
2h_{_\C}(v,w):=h(v,\overline{w})+\overline{h(\overline{v},w)},\quad v,w\in V\otimes_{_\R}\C,
$$
where we still note by  $h$ the $\C$-linear extension of $h$.        
Consider now the complex vector bundles $F:=\Lambda ^{1,0}_{_{J}}T^*_{_X},\, E:=F^{\otimes 2}\otimes \overline{F}$  and the hermitian vector bundles
\begin{eqnarray*}
(E,\left<\cdot,\cdot\right>_{\omega})&:=&(F,h^*)\otimes(F,h^*)\otimes(\overline{F},h^*)
\\
(E,\left<\cdot,\cdot\right>_{\varphi})&:=&(F,h^*_\varphi)\otimes(F,h^*_\varphi)\otimes(\overline{F},h^*_\varphi)
\\
(E,\left<\cdot,\cdot\right>_{\omega,\varphi})&:=&(F,h^*)\otimes(F,h^*_\varphi)\otimes(\overline{F},h^*_\varphi)
\\
(E,\left<\cdot,\cdot\right>_{\varphi,\omega})&:=&(F,h^*_\varphi)\otimes(F,h^*_\varphi)\otimes(\overline{F},h^*).
\end{eqnarray*}
For example the last two hermitian metrics are expressed in local coordinates by the expressions
\begin{eqnarray*}
\left<\alpha,\beta\right>_{\omega,\varphi}
&=&
2^{-3}h^*(dz_p,d\bar{z}_q)\,h^*_{\varphi}(dz_j,d\bar{z}_l)\,\overline{h^*_{\varphi}(dz_k,d\bar{z}_m)}\,
\alpha_{pj\bar{k}}\overline{\beta_{ql\bar{m}}}
\\
&=& 2^3\omega^{q\bar{p}}\omega^{l\bar{j}}_{\varphi}\omega^{k\bar{m}}_{\varphi}\alpha_{pj\bar{k}}\overline{\beta_{ql\bar{m}}},
\\
\left<\alpha,\beta\right>_{\varphi,\omega}
&=& 2^3\omega^{q\bar{p}}_{\varphi}\omega^{l\bar{j}}_{\varphi}\omega^{k\bar{m}}\alpha_{pj\bar{k}}\overline{\beta_{ql\bar{m}}}
\end{eqnarray*}
where $\alpha=\alpha_{pj\bar{k}}dz_p\otimes dz_j\otimes d\bar{z}_k$ and $\beta=\beta_{pj\bar{k}}dz_p\otimes dz_j\otimes d\bar{z}_k$. With such notations we can state the following lemma (see also \cite{Yau}) that we will prove at the end of the section.
\begin{lem}\label{C3CalaBox}
Let $(X,\omega )$ be a polarized Fano manifold of complex dimension $n$ with $\omega\in 2\pi c_1$ and let $(\omega_t)_t$ be the K\"ahler-Ricci flow. Suppose that there exist constants $k,\, K>0$ such that 
$k^{-1}\omega\leq \omega _t\leq K\omega$
for all times $t\geq 0$.
Then there exist constants $C_1,\,C_2>0$ depending only on the constants $k$, $K$ and $\omega$, such that the uniform estimate 
\begin{eqnarray*}
\Openbox_{\,t}|\nabla^{1,0}_{\omega}\partial\bar{\partial}\varphi_t|^2_t
\geq
-C_1
|\nabla^{1,0}_{\omega}\partial\bar{\partial}\varphi_t|^2_t -C_2
\end{eqnarray*}
holds for every $t\geq 0$.
\end{lem}  
By using Perelman's uniform estimate $|\dot{\varphi}_t|\leq C$ and a slight modification of Yau's computation of the $C^2$ and $C^3$ uniform estimates for the complex Monge-Amp\`ere operator (see also \cite{Cao}), we find the following result.
\begin{prop}\label{CalabiC3KRicc}
Let $X$ be a Fano manifold of complex dimension $n\geq 2$ and let $\omega \in 2\pi c_1(X)$ be a K\"ahler metric. If the K\"ahler-Ricci flow
$$
\dot{\varphi}_t=\log\,\frac{\omega _t ^n}{\omega ^n}+\varphi_t+c_t-h_{\omega},  
$$
%with initial value $\varphi_0=0$ 
satisfies the uniform estimate $\omega^n_t\geq K_0\,\omega^n$ 
for some constant $K_0>0$ independent of $t\in [0,+\infty)$,
then there exist positive constants $k_0,\,K,\, K'>0$ independents of $t\in [0,+\infty)$, such that  the uniform estimates 
$0<2n+\Delta _{\omega}\varphi_t \leq K$, $|\partial\bar{\partial}\varphi_t|_{\omega} < (K+2\sqrt{n})/2$, 
$k^{-1}_0\omega <\omega _t<(K/2)\omega$ and $|\nabla^{1,0}_{\omega}\partial\bar{\partial}\varphi_t|_{\omega}\leq K'$ holds for all $t\in [0,+\infty)$.
\end{prop}  
(The $C^2$-uniform estimate is obvious in the case $n=1$.)
\\
\\
$Proof$. We define the operator $\Openbox_{\,t}:=\Delta _t-2\frac{\partial}{\partial t}$. Consider the smooth function $A:=\log(2n+\Delta _{\omega}\varphi_t)-k(\varphi_t+c_t)$ over $X\times [0,+\infty)$, where the constant $k$ will be chosed later. We have the equality
\begin{eqnarray}\label{BoxA} 
\Openbox_{\,t} A= \frac{ \Openbox_{\,t}\,\Delta _{\omega}\varphi_t }{2n +\Delta_{\omega }\varphi_t}
-
\frac{2|\partial \Delta _{\omega}\varphi_t  |^2_t}{(2n +\Delta_{\omega }\varphi_t)^2}
-
k\,\Openbox_{\,t}\varphi_t+2ka_t.  
\end{eqnarray}
Set $C:=\min_{x\in X}\lambda ^{\omega}_1(x)$. Using the proposition \ref{RiccCaYau} and the fact that $\Tr_{\varphi}\omega>0$ we find the inequality
\begin{eqnarray*} 
2\Tr_{\omega} \Ric(\omega_{\varphi})\geq  
-\Delta _{\varphi}\Delta _{\omega}\varphi 
+4C(2n +\Delta_{\omega }\varphi )\Tr_{\varphi}\omega
+\frac{2|\partial \Delta _{\omega}\varphi  |^2_{\varphi }}{2n +\Delta_{\omega }\varphi}
\end{eqnarray*}  
which combined with the equality \eqref{BoxA} gives
\begin{eqnarray}\label{BoxA1} 
\Openbox_{\,t} A\geq  -2\,\frac{\Delta_{\omega }\dot{\varphi}_t+ \Tr_{\omega} \Ric(\omega_t)}{2n +\Delta_{\omega }\varphi_t}
+
4C\Tr_t\omega
-
k\,\Openbox_{\,t}\varphi_t+2ka_t,  
\end{eqnarray}
where $\Tr_t\omega$ is the trace of $\omega$ respect to the metric $\omega_t$. Taking the trace of the K\"ahler-Ricci flow identity $i\partial\bar{\partial}\dot{\varphi}_t=\omega_t-\Ric(\omega_t)$
respect to $\omega$, we find the equality
\begin{eqnarray}\label{DelKRicci}
\Delta _{\omega}\dot{\varphi}_t =2n+\Delta _{\omega}\varphi_t-\Tr_{\omega} \Ric(\omega_t).
\end{eqnarray}
Moreover concerning the term $\Openbox_{\,t}\varphi_t$, we remark the trivial identity $\Delta _t\varphi_t =-\Tr_t\omega+2n$. Using this identity with the equality \eqref{DelKRicci} in the inequality \eqref{BoxA1}, we find
\begin{eqnarray}\label{DelA2} 
\Openbox_{\,t}A\geq  -2+(4C+k)\Tr_t\omega+2k(\dot{\varphi}_t-n+a_t).
\end{eqnarray}
Consider now the trivial inequality $\sum_{l=1}^n b_1\dots\widehat{b}_l\dots b_n\leq (\sum_{l=1}^n b_l)^{n-1}$ for any positive number $b_l$. Taking the $1/(n-1)$-th power of this inequality with the terms
\\
$b_l:=1/(1+2\partial^2_{l\bar{l}}\,\varphi_t)$ we find, in $\omega$-orthogonal and $\omega_t$-diagonal coordinates in a point $x$, the expressions
\begin{eqnarray*} 
\frac{\Tr_t\omega}{4}=\sum_l\frac{1}{1+2\partial^2_{l\bar{l}}\,\varphi_t}
&\geq& 
\left(\frac
{\sum_l\left(1+2\partial^2_{l\bar{l}}\,\varphi_t\right)}
{\prod_l \left(1+2\partial^2_{l\bar{l}}\,\varphi_t\right)}  \right)^{\frac{1}{n-1}}
\\
&=&
K_n\,e^{\frac{\varphi_t+c_t -h-\dot{\varphi}_t}{n-1}}(2n+\Delta _{\omega}\varphi )^{\frac{1}{n-1}},     
\end{eqnarray*}
where $K_n:=2^{\frac{-1}{n-1}}>0$. We choose $k$ such that $(4C+k)=4^{-1}$ and we consider the function $u:=e^A=(2n+\Delta _{\omega}\varphi_t )e^{-k(\varphi_t+c_t)}$. Then the previous inequality combined with the Perelman uniform estimate $|\dot{\varphi}_t|\leq C'$ and with the estimate $|\varphi_t+c_t|\leq C$, gives
$$
\frac{\Tr_t\omega}{4}\geq K_n e^{\frac{(1+k)(\varphi_t+c_t) -h-\dot{\varphi}_t}{n-1}}u^{\frac{1}{n-1}}\geq C_0\, u^{\frac{1}{n-1}}
$$ 
for some constant $C_0>0$ independent of $t$.  Then the inequality \eqref{DelA2} reduces to the inequality 
\begin{eqnarray}\label{UnifKRest} 
\Openbox_{\,t} A\geq -C_1+C_0\,u^{\frac{1}{n-1}},  
\end{eqnarray}
with $C_0,\,C_1>0$. For all $T>0$, a point  $(x_0,t_0)$ is a maximum point for $A$ in $X\times[0,T]$ if and only if is also a maximum point for $u$ in $X\times[0,T]$. If $t_0=0$ then $u\leq C_2$ over $X\times[0,T]$, with $C_2>0$ independent of $T$. 
If not $\frac{\partial A}{\partial t}(x_0,t_0)\geq 0$ and $\Delta _tA(x_0,t_0)\leq 0$. Using the inequality \eqref{UnifKRest}  we find 
$u(x_0,t_0)\leq C_3$, where the constant $C_3>0$ is independent of $T$. This implies the estimate
$u\leq \max\{C_2,C_3\}$ on $X\times[0,+\infty)$. So in conclusion we have found the required a priori estimate 
$0<2n+\Delta _{\omega }\varphi_t\leq K$. Moreover $2|\omega_t |_{\omega}<\Tr_{\omega}\omega_t$, since $\omega_t>0$. This
implies the required a priori estimate $|\partial\bar{\partial}\varphi_t|_{\omega} < (K+2\sqrt{n})/2$. The inequality
$0<2+4\partial^2_{l\bar{l}}\,\varphi_t<2n+\Delta _{\omega }\varphi_t\leq K$ implies $\omega_t<(K/2)\omega$. 
By using the hypothesis we find
$$
K_0\leq \omega _t ^n/\omega ^n=\prod_l \left(1+2\partial^2_{l\bar{l}}\,\varphi_t\right)<(K/2)^{n-1}\left(1+2\partial^2_{s\bar{s}}\varphi_t\right),
$$
%By definition of the K\"ahler-Ricci flow we have 
%$$
%K_0\leq e^{\dot{\varphi}_t+h-(\varphi_t+c_t)}=\omega _t ^n/\omega ^n=\prod_l \left(1+2\partial^2_{l\bar{l}}\,\varphi_t\right)<(K/2)^{n-1}\left(1+2\partial^2_{s\bar{s}}\varphi_t\right),
%$$
for all $s$, which implies $k^{-1}_0\omega< \omega _t$ for some uniform constant $k_0>0$.
Then by lemma \ref{C3CalaBox} we deduce the estimate
\begin{eqnarray}\label{C3CalKRF}
\Openbox_{\,t}|\nabla^{1,0}_{\omega}\partial\bar{\partial}\varphi_t|^2_t
\geq
-C_1
|\nabla^{1,0}_{\omega}\partial\bar{\partial}\varphi_t|^2_t -C_2\,.
\end{eqnarray}
The equality \eqref{TrRic} proved in the proposition \ref{RiccCaYau} gives the intrinsic identity
\begin{eqnarray*}
2\Tr_{\omega} \Ric(\omega_{\varphi})= 
-\Delta _{\varphi}\Delta _{\omega}\varphi 
+2\Tr_{\varphi}(\omega_{\varphi}\cdot\Rm_{\omega})
+4|\nabla^{1,0}_{\omega}\partial\bar{\partial}\varphi|^2_{\omega,\varphi},
\end{eqnarray*}
where $%-4\omega_{\varphi}(\omega^*_{\varphi}\contract \Rm_{\omega})
2\Tr_{\varphi}(\omega_{\varphi}\cdot\Rm_{\omega})
\geq 4\lambda^{\omega}_1(2n +\Delta_{\omega }\varphi )\Tr_{\varphi}\omega\geq -C_3$, with $C_3>0$. So using the identity \eqref{DelKRicci}, we deduce the inequality
$$
\Openbox_{\,t}\Delta _{\omega}\varphi_t\geq(4/k_0)|\nabla^{1,0}_{\omega}\partial\bar{\partial}\varphi_t|^2_t-C_4,
$$
$C_4>0$. By taking $C_5:=k_0(C_1+1)/4>0$ we deduce by \eqref{C3CalKRF} the estimate
\begin{eqnarray}\label{FinalC3Est}
\Openbox_{\,t}\left(|\nabla^{1,0}_{\omega}\partial\bar{\partial}\varphi_t|^2_t+C_5\Delta _{\omega}\varphi_t\right)\geq |\nabla^{1,0}_{\omega}\partial\bar{\partial}\varphi_t|^2_t-C_6,
\end{eqnarray}
$C_6>0$. For all $T>0$ consider a maximum point $(x_0,t_0)$ for the function $B:=|\nabla^{1,0}_{\omega}\partial\bar{\partial}\varphi_t|^2_t+C_5\Delta _{\omega}\varphi_t$ over $X\times [0,T]$. As before we can assume $t_0>0$, which implies $\Openbox_{\,t} B(x_0,t_0)\leq 0$. Then the estimate \eqref{FinalC3Est} implies the inequality 
$B(x_0,t_0)\leq C_6+C_5\Delta _{\omega}\varphi_{t_0}(x_0)\leq C_7$, for some constant $C_7>0$ independent of $T$. We deduce in conclusion the required third order uniform estimate
$|\nabla^{1,0}_{\omega}\partial\bar{\partial}\varphi_t|^2_{\omega}\leq K'$, since the metrics $\omega_t$ are uniformly equivalent to the initial metric $\omega$.
\hfill$\Box$
\\
\\
{\bf Proof of lemma \ref {C3CalaBox}}. In order to avoid confusion with notations in the computations that will follow we will note $\varphi:=\varphi_t$ and $\omega_{\varphi}:=\omega_t$. This will apply until equality \eqref{C3Expres}.
According to lemma \ref{KahlCurvGeod3} let consider $(z_1,...,z_n)$ be $\omega$-geodesic holomorphic coordinates of third order, with center a point $x$ such that the metric $\omega _{\varphi}$ can be written in diagonal form in $x$. Explicitly 
$\omega =\frac{i}{2}\omega _{l\bar{r}}\,dz_l\wedge d\bar{z}_r$, where 
\begin{eqnarray*}
&\displaystyle{
\omega _{l\bar{r}}=\delta _{lr}-C^{j\bar{k}}_{rl}z_j\bar{z}_k-C^{pj\bar{k}}_{l\bar{r}}z_pz_j\bar{z}_k-\overline{C^{pj\bar{k}}_{r\bar{l}}}z_k\bar{z}_p\bar{z}_j+O(|z|^4), }&
\\
&\displaystyle{
2R_{j\bar{k}l\bar{r}}(x) =
 C^{j\bar{k}}_{rl},
 \quad
2\nabla^{1,0}_pR_{j\bar{k}l\bar{r}}(x)=C^{pj\bar{k}}_{l\bar{r}},
\quad
2\nabla^{0,1}_{\bar{p}}R_{j\bar{k}l\bar{r}}(x)=\overline{C^{pj\bar{k}}_{r\bar{l}}},}&
\\
&\displaystyle{
\overline{C^{j\bar{k}}_{rl}}=C^{k\bar{j}}_{lr},\quad C^{j\bar{k}}_{rl}=C^{l\bar{k}}_{rj}=C^{j\bar{r}}_{kl},}&
\end{eqnarray*}
and the coefficients $C^{p,j,\bar{k}}_{l\bar{r}}$ are symmetric respect to the indexes $p,j,l$ and $k,r$.
We define $a_l:=\omega^{l\bar{l}}_{\varphi}$.
By deriving the Ricci tensor $\Ric(\omega_{\varphi})=\Ric(\omega_{\varphi})_{j\bar{k}}\,dz_j\wedge d\bar{z}_k$,
\begin{eqnarray*}
\Ric(\omega_{\varphi})_{j\bar{k}}
=
-i
(
\partial^2_{j\bar{k}}\omega_{p\bar{r}} +2\varphi_{jp \bar{k}\bar{r}}  
)\omega^{r\bar{p}}_{\varphi}
+
i(\partial_j\omega_{p\bar{s}}+2\varphi_{jp\bar{s}})\omega^{s\bar{t}}_{\varphi}
(\partial_{\bar{k}}\omega_{t\bar{r}}+2\varphi_{t\bar{r}\bar{k}})\omega^{r\bar{p}}_{\varphi},
\end{eqnarray*}
we find at the point $x$ the expression
\begin{eqnarray}\label{CovDerRic}
\nabla^{1,0}_{\omega,t}\Ric(\omega_{\varphi})_{j\bar{k}}
&=&
ia_p
%\left(
( \underbrace{C^{tj\bar{k}}_{p\bar{p}}}_{C1}
\, - \, 
\underbrace{2\varphi_{jtp\bar{p}\bar{k}}}_{A1} 
%\right)
)
\, + \, \underbrace{4ia_pa_l\varphi_{tjp\bar{l}} 
\varphi_{l\bar{p}\bar{k}}}_{A2}  \nonumber
\\\nonumber
&+&
2ia_pa_l
%\left[
\Big[
%\left(
(
\underbrace{C^{j\bar{k}}_{lp}}_{C2} \,
+ \, \underbrace{2\varphi_{jp\bar{k}\bar{l}}}_{A3} 
%\right)
)
\varphi_{tl\bar{p}} 
\, + \, 
%\left(
(
\underbrace{2\varphi_{tl\bar{p}\bar{k}}}_{A4} 
\, - \, \underbrace{C^{t\bar{k}}_{pl}}_{C3} 
%\right)
)
\varphi_{jp\bar{l}} 
%\right]
\Big]\nonumber
\\
&-&
8ia_pa_la_r
%\left(
(
\underbrace{\varphi_{jp\bar{r}}\varphi_{tr\bar{l}}\varphi_{l\bar{p}\bar{k}}}_{A5}   
\, + \,
\underbrace{\varphi_{jp\bar{l}}\varphi_{tr\bar{p}}\varphi_{l\bar{r}\bar{k}}}_{A6}   
%\right)
).
\end{eqnarray}
The utility of the underbraces will be discussed later. Consider now the tensor $\nabla^{1,0}_\omega\partial \bar{\partial} \varphi=\alpha_{p k \bar{l}}\,dz_p\otimes(dz_k\wedge d\bar{z}_l)$, where
$$
\alpha_{p k \bar{l}}  : =   \varphi_{p k \bar{l}}  -  
\partial_p  \, \omega_{k \bar{r}}  \, \omega^{r \bar{s}} \varphi_{s
\bar{l}}
$$
and the derivative of its norm
\begin{eqnarray*}
\partial_p 
| \nabla^{1,0}_\omega 
\partial \bar{\partial} \varphi|^2_\varphi
& = & 
2^3 \partial_p \, 
\omega^{s \bar{t}}_\varphi  
\omega^{l \bar{j}}_\varphi \,
\omega^{k \bar{r}}_\varphi \, 
\alpha_{t j \bar{k}} \,
\overline{\alpha_{s l \bar{r}}}  
+ 
2^3 \, \omega^{s \bar{t}}_\varphi \partial_p \, \omega^{l
\bar{j}}_\varphi \, \omega^{k \bar{r}}_\varphi \, \alpha_{t j \bar{k}}
\, \overline{\alpha_{s l \bar{r}}} 
\\
& + &   2^3 \, \omega^{s \bar{t}}_\varphi \,
\omega^{l \bar{j}}_\varphi \, \partial_p \, \omega^{k
\bar{r}}_\varphi \, 
\alpha_{t j \bar{k}} \, 
\overline{\alpha_{s l \bar{r}}} 
+ 
2^3 \, \omega^{s \bar{t}}_\varphi \, 
\omega^{l \bar{j}}_\varphi
\,
\omega^{k \bar{r}}_\varphi \, 
\partial_p  \alpha_{t j \bar{k}}  \,
\overline{\alpha_{s l \bar{r}}} 
\\
& + & 2^3 \, \omega^{s t}_\varphi \, 
\omega^{l \bar{j}}_\varphi \,
\omega^{k \bar{r}}_\varphi \, 
\alpha_{t j \bar{k}} \, \partial_p \, 
\overline{\alpha_{s l \bar{r}}}  \,  
\end{eqnarray*}
By using the expressions of the derivatives
\begin{eqnarray*}
\partial_p  \alpha_{t j \bar{k}} 
&=& 
 \varphi_{p t j \bar{k}} \, -
\partial^2_{t p} \,
\omega_{j \bar{a}} \, \omega^{a \bar{b}} \,
\varphi_{b \bar{k}} 
-
 \partial_t \,
\omega_{j \bar{a}} \,
\partial_p \, \omega^{a \bar{b}} \, \varphi_{b \bar{k}} 
- \partial_t \, \omega_{j \bar{a}} \, \omega^{a \bar{b}} \,
\varphi_{b p \bar{k}} 
\\
\partial_p \overline{\alpha_{s l \bar{r}}} 
&=&
\varphi_{p r \bar{s} \bar{l}} \, - 
\partial^2_{p \bar{s}} \, 
\omega_{a  \bar{l}} \,
\omega^{b \bar{a}} \,
\varphi_{r \bar{b}} 
-
\partial_{\bar{s}} \,
\omega_{a \bar{l}} \, \partial_p \,
\omega^{b \bar{a}} \,
\varphi_{r \bar{b}} 
-
\partial_{\bar{s}} \,
\omega_{a \bar{l}} \,
\omega^{b , \bar{a}} \,
\varphi_{p r \bar{b}},
\end{eqnarray*}
we find the the following expression for the Laplacian at the point $x$.
\begin{eqnarray*}
\Delta_\varphi |
\nabla^{1,0}_\omega 
\partial \bar{\partial} 
\varphi |^2_\varphi & = & 
2^5 \, 
\omega^{q \bar{p}}_\varphi \,
\left[
\partial^2_{p \bar{q}} \,
\omega^{s \bar{t}}_\varphi \,
\omega^{l \bar{j}}_\varphi \,
\omega^{k r}_\varphi \,
\varphi_{t j \bar{k}} 
\varphi_{r \bar{s} \bar{l}} \right. \\
& + & \partial_p \,
\omega^{s \bar{t}}_\varphi \,
\partial_{\bar{q}} \,
\omega^{l \bar{j}}_\varphi \,
\omega^{k \bar{r}}_\varphi \,
\varphi_{t j \bar{k}} \,
\varphi_{r \bar{s} \bar{l}} \, + \,
\partial_p \, 
\omega^{s \bar{t}}_\varphi \,
\omega^{l \bar{j}}_\varphi \,
\partial_{\bar{q}} \,
\omega^{k \bar{r}}_\varphi \, 
\varphi_{t j \bar{k}} \,
\varphi_{r \bar{s} \bar{l}} \\
& + & \partial_p \, 
\omega^{s , \bar{t}}_\varphi \,
\omega^{l \bar{j}} _\varphi \,
\omega^{k \bar{r}}_\varphi ( \varphi_{t j \bar{q} \bar{k}} \, + \, C^{t
\bar{q}}_{a j} \, \varphi_{a \bar{k}} ) \, \varphi_{r \bar{s}
\bar{l}} \\
& + & \partial_p \, 
\omega^{s \bar{t}}_\varphi \,
\omega^{l \bar{j}}_\varphi \,
\omega^{k \bar{r}}_\varphi \,
\varphi_{t j \bar{k}} \,
\varphi_{r \bar{s}\bar{l} \bar{p}} \, + \,
\partial_{\bar{q}} \,
\omega^{s \bar{t}}_\varphi \,
\partial_p \, 
\omega^{l \bar{j}}_\varphi \,
\omega^{l \bar{r}}_ \varphi \,
\varphi_{t j \bar{k}} \, 
\varphi_{r \bar{s} \bar{l}} \\
& + & \omega^{s \bar{t}}_\varphi \, 
\partial^2_{p \bar{q}} \,
\omega^{l \bar{j}}_\varphi \, 
\omega^{k \bar{r}}_\varphi \, \varphi_{t j 
\bar{k}} \,
\varphi_{r \bar{t} \bar{l}} \, + \,
\omega^{s \bar{t}}_\varphi \, \partial_p \,
\omega^{l \bar{j}}_\varphi \, \partial_{\bar{q}} \,
\omega^{k \bar{r}}_\varphi \, \varphi_{t j \bar{k}} \, \varphi_{r
\bar{s} \bar{l}} 
\end{eqnarray*}
\begin{eqnarray*}
& + & \omega^{s \bar{t}}_\varphi \, 
\partial_p \,
\omega^{l \bar{j}}_\varphi \,
\omega^{k \bar{r}}_\varphi \,
( \varphi_{t j \bar{q} \bar{k}} \, + \,
C^{t \bar{q}}_{s j} \, 
\varphi_{s \bar{k}} ) \, \varphi_{r \bar{s}{l}} \\
& + & 
\omega^{s \bar{t}}_\varphi \,
\partial_p \, \omega^{l \bar{j}}_\varphi \,
\omega^{k \bar{r}}_\varphi \,
\varphi_{t j \bar{k}} \,
\varphi_{r \bar{s} \bar{l} \bar{q}} \, + \,
\partial_{\bar{q}} \,
\omega^{s \bar{t}}_\varphi \,
\omega^{l \bar{j}}_\varphi \,
\partial_p \, \omega^{k \bar{r}}_\varphi \, \varphi_{t j \bar{k}} \,
\varphi_{r \bar{s} \bar{l}} \\
& + & \omega^{s \bar{t}}_\varphi \,
\partial_{\bar{q}} \,
\omega^{l \bar{j}}_\varphi \,
\partial_p \,
\omega^{k \bar{r}}_\varphi \,
\varphi_{t j \bar{k}} \,
\varphi_{r \bar{s} \bar{l}} \, + \,
\omega^{s \bar{t}}_\varphi \,
\omega^{l \bar{j}}_\varphi \, \partial^2_{p \bar{q}} \, \omega^{k
\bar{r}}_\varphi \, \varphi_{t j \bar{k}} \,
\varphi_{r \bar{s} \bar{l}} \\
& + & \omega^{s \bar{t}}_\varphi
\omega^{l \bar{j}}_\varphi
\partial_p \,
\omega^{k \bar{r}}_\varphi 
( \varphi_{t j \bar{k} \bar{q}}  + 
C^{t \bar{q}}_{s j} 
\varphi_{s \bar{k}} )
\varphi_{r \bar{s}\bar{l}} + 
\omega^{s \bar{t}}_\varphi
\omega^{l \bar{j}}_\varphi 
\partial_p \,
\omega^{k \bar{r}}_\varphi 
\varphi_{t j \bar{k}} 
\varphi_{r \bar{s} \bar{l} \bar{q}} \\
& + & \partial_{\bar{q}} \,
\omega^{s , \bar{t}}_\varphi \,
\omega^{l \bar{j}}_\varphi \,
\omega^{k \bar{r}}_\varphi\,
\varphi_{p t j \bar{k}} \,
\varphi_{r \bar{s} \bar{l}} \, + \,
\omega^{s , \bar{t}}_\varphi \,
\partial_{\bar{q}} \,
\omega^{l \bar{j}}_\varphi \,
\omega^{k \bar{r}}_\varphi \,
\varphi_{p t j \bar{k}} \, \varphi_{r \bar{s} \bar{l}} \\
& + & \omega^{s \bar{t}}_\varphi \,
\omega^{l \bar{j}}_\varphi \, \partial_{\bar{q}} \, \omega^{k
\bar{r}}_\varphi \, 
\varphi_{p t j \bar{k}} \, \varphi_{r \bar{s} \bar{l}}  \\ 
& + & 
\omega^{s \bar{t}}_ \varphi \,
\omega^{l j }_\varphi \,
\omega^{k \bar{r}}_\varphi ( \varphi_{p t j \bar{k} \bar{q}} \, + \,
C^{t p \bar{q}}_{j {\bar{a}}} \,
\varphi_{a \bar{k}} \, + \,
C^{t \bar{q}}_{a j} \, 
\varphi_{a p \bar{k}} ) \, 
\varphi_{r \bar{s} \bar{l}} \\
& + &
\omega^{s  \bar{t}}_\varphi \,
\omega^{l \bar{j}}_\varphi \,
\omega^{k \bar{r}}_\varphi 
\varphi_{p t j \bar{k}} \,
\varphi_{r \bar{s} \bar{l} \bar{q}} \, + \,
\partial_{\bar{q}} \,
\omega^{s  \bar{t}}_\varphi \,
\omega^{l \bar{j}}_\varphi \,
\omega^{k \bar{r}}_\varphi \varphi_{t j \bar{k}} ( \varphi_{p r
\bar{s} \bar{l}} +  C^{p \bar{s}}_{l a} 
\varphi_{r \bar{a}} ) \\
& + & \omega^{s \bar{t}}_\varphi \,
\partial_{\bar{q}} \, \omega^{l \bar{j}}_\varphi \, \omega^{k
\bar{r}}_\varphi \, 
\varphi_{t j \bar{k}} \, ( \varphi_{p r \bar{s} \bar{l}} \, + \, C^{p
\bar{s}}_{l a} \, \varphi_{r \bar{a}} ) \\
& + & 
\omega^{s  \bar{t}}_\varphi \,
\omega^{l \bar{j}}_\varphi \, \partial_{\bar{q}} \,
\omega^{k \bar{r}}_\varphi \,
\varphi_{t j \bar{k}} ( \varphi_{p r \bar{s} \bar{l}} \, + \, C^{p
\bar{s}}_{l a} \, \varphi_{r \bar{a}} ) \\
& + & \omega^{s  \bar{t}}_\varphi \,
\omega^{l \bar{j}}_\varphi \,
\omega^{k \bar{r}}_\varphi \,
( \varphi_{t j \bar{k} \bar{q}} \, + \,
C^{t \bar{q}}_{a j} \,
\varphi_{a \bar{k}} ) \,
( 
\varphi_{p r \bar{s} \bar{l}} \, + \, C^{p \bar{s}}_{l a} \,
\varphi_{r \bar{a}} ) \\
& + & \omega^{s  \bar{t}}_\varphi \,
\omega^{l j}_{\varphi} \,
\omega^{k \bar{r}}_\varphi \, \varphi_{t j \bar{k}} \,
\left(
\varphi_{p r \bar{s} \bar{l} \bar{q}} \, + \,
\overline{C^{t q \bar{p}}_{l \bar{a}} } \, 
\varphi_{r \bar{a}} \, + \, 
C^{p \bar{t}}_{l a} \, 
\varphi_{r \bar{a} \bar{q}} 
\right)
\Big].
\end{eqnarray*}
Then using the expressions
$
\partial_{\bar{l}} \,
\omega^{s \bar{t}}_\varphi =
- 2 a_s a_t \, 
\varphi_{s \bar{l \bar{t}}} 
$
and
$$
\partial^2_{k \bar{l}} \,
\omega^{s  \bar{t}}_\varphi  = 
a_s a_t 
\left[
C^{k \bar{l}}_{t s} - 2 
\varphi_{k s \bar{l} \bar{t}} \, + \,
4 a_r 
\left(
\varphi_{k s \bar{r}} 
\varphi_{r \bar{l} 
\bar{t}} \, + \,
\varphi_{k r \bar{t}} \, 
\varphi_{s \bar{l \bar{r}}} 
\right) 
\right],
$$
at the point $x$ of the derivatives of the inverse matrixs
we find the following expression for the Laplacian at the point $x$.
\begin{eqnarray*}
\Delta _{\varphi}|\nabla^{1,0}_{\omega}\partial\bar{\partial}\varphi|^2_{\varphi }
&=&
2^5a_pa_ta_ja_k
\Big\{
\\
&&
a_l
%\left[
[4a_r
%\left(
(
\varphi_{pl\bar{r}}  \varphi_{ r\bar{p}\bar{t}} \,  + \, 
\underbrace{\varphi_{pr\bar{t}}  \varphi_{ l\bar{p}\bar{r}}}_{B1d} 
%\right)
)
- \, 2\varphi_{pl\bar{p}\bar{t}}  +C^{p\bar{p}}_{tl}
%\right]
]
\varphi_{tj\bar{k}}
\varphi_{ k\bar{j}\bar{l}}   
\\
&+&
4a_la_r
%\left(
(
\underbrace{\varphi_{k\bar{r}\bar{l}} \varphi_{ l\bar{p}\bar{j}}}_{A5} 
\, + \, \underbrace{\varphi_{ l\bar{r}\bar{j}} \varphi_{
k\bar{p}\bar{l}}}_{B1d} 
%\right)
)
\varphi_{pr\bar{t}}\varphi_{tj\bar{k}}
\\
&-&
2a_l
%\left(
(
\underbrace{\varphi_{tj\bar{p}\bar{k}}}_{\overline{A4}} 
\, + \,C^{t\bar{p}}_{kj}\varphi_{k\bar{k}}  
%\right)
)
\varphi_{pl\bar{t}}\varphi_{k\bar{l}\bar{j}} \, - \,
\underbrace{2a_l \varphi_{k\bar{l}\bar{j}\bar{p}}\varphi_{pl\bar{t}}
\varphi_{tj\bar{k}}}_{B1b}
\\
&+&
\underbrace{4a_la_r\varphi_{pl\bar{j}} \varphi_{tj\bar{k}} 
\varphi_{ k\bar{r}\bar{l}} \varphi_{ r\bar{p}\bar{t}}}_{B1c} 
\\
&+&
a_l
%\left[
[
4a_r
%\left(
( \varphi_{pl\bar{r}} 
 \varphi_{r\bar{p}\bar{j}} \, + \,
\underbrace{\varphi_{pr\bar{j}}\varphi_{l\bar{p}\bar{r}}}_{B1c}
%\right)
)
\, - \, 2\varphi_{ pl\bar{p}\bar{j}}+C^{p\bar{p}}_{jl} 
%\right]
]
\varphi_{tj\bar{k}}
\varphi_{k \bar{t}\bar{l}}
%}_{Hey}
\\
&+&
4a_la_r
%\left(
(
\underbrace{\varphi_{pl\bar{j}}\varphi_{k\bar{p}\bar{r}}}_{A6} \, + \,
\underbrace{\varphi_{pk\bar{r}}\varphi_{l\bar{p}\bar{j}}}_{\overline{A5}} 
%\right)
)
\varphi_{tj\bar{k}}\varphi_{r\bar{t}\bar{l}}
\\
&-&
2a_l
%\left(
( 
\underbrace{\varphi_{tj\bar{p}\bar{k}}}_{A4} \, + \, 
C^{t\bar{p}}_{kj}\varphi_{k\bar{k}}  
%\right) 
)
\varphi_{pl\bar{j}}\varphi_{k\bar{t}\bar{l}}
\, - \,
\underbrace{2a_l\varphi_{k\bar{t}\bar{l}\bar{p}}
\varphi_{pl\bar{j}}\varphi_{tj\bar{k}}}_{B1a} 
\\
&+&
\underbrace{4a_la_r\varphi_{pk\bar{r}}
\varphi_{tj\bar{k}}\varphi_{r\bar{l}\bar{t}}
\varphi_{l\bar{p}\bar{j}} }_{\overline{A6}} 
\\
&+&
a_l
%\left[
[
4a_r
%\left(
(
\varphi_{pk\bar{r}} \varphi_{r\bar{p}\bar{l}}
\, + \, \underbrace{ \varphi_{pr\bar{l}}  \varphi_{k
\bar{p}\bar{r}}}_{B2a}   
%\right)
)
-2\varphi_{pk\bar{p}\bar{l}}  +C^{p\bar{p}}_{lk}
%\right]
]
\varphi_{tj\bar{k}} \varphi_{l\bar{t}\bar{j}}
\end{eqnarray*}
\begin{eqnarray*}
&-&
\underbrace{2a_l
%\left( 
(
\varphi_{tj\bar{p}\bar{k}}+C^{t\bar{p}}_{kj} \varphi_{k\bar{k}}
%\right)
)
\varphi_{pk\bar{l}} \varphi_{l\bar{t}\bar{j}} }_{B2}
\, - \, 
\underbrace{2a_l\varphi_{l\bar{t}\bar{j}\bar{p}}
\varphi_{pk\bar{l}}\varphi_{tj\bar{k}}}_{\overline{A2}} 
\\
&-&
\underbrace{2a_l\varphi_{ptj\bar{k}}\varphi_{l\bar{p}\bar{t}}\varphi_{k\bar{j}\bar{l}}}_{B1} 
\, - \,
\underbrace{2a_l\varphi_{ptj\bar{k}}
\varphi_{l\bar{p}\bar{j}}\varphi_{k\bar{t}\bar{l}} }_{A2} 
\, - \, 
\underbrace{2a_l\varphi_{ptj\bar{k}}
\varphi_{k\bar{p}\bar{l}} \varphi_{l\bar{t}\bar{j}}}_{B1} 
\\
&+&
%\left(
( \underbrace{\varphi_{ptj\bar{k}\bar{p}}}_{A1} 
\, + \, 
C^{t\bar{p}}_{rj} \varphi_{rp\bar{k}} \,  + \, 
\underbrace{C^{tp\bar{p}}_{j\bar{k}} \varphi_{k\bar{k}} }_{C1} 
%\right)
)
\varphi_{k\bar{t}\bar{j}}  
\, + \, 
\underbrace{\varphi_{ptj\bar{k}}\varphi_{k\bar{p}\bar{t}\bar{j}}}_{B1}
\\
&-&
2a_l
%\left(
(
\underbrace{\varphi_{pk\bar{l}\bar{j}}}_{\overline{A3}} 
\, + \, C^{p,\bar{l}}_{jk}\varphi_{k\bar{k}}
%\right)
)
\varphi_{tj\bar{k}}\varphi_{l\bar{p}\bar{t}} 
\\
&-&
2a_l
%\left(
(
\underbrace{\varphi_{pk\bar{t}\bar{l}}}_{A3} 
\, + \, C^{p\bar{t}}_{lk}\varphi_{k\bar{k}} 
%\right)
)
\varphi_{tj\bar{k}}\varphi_{l\bar{p}\bar{j}} 
\\
&-&
\underbrace{
2a_l
%\left(
(
\varphi_{pl\bar{t}\bar{j}} +C^{p\bar{t}}_{jl}\varphi_{l\bar{l}} 
%\right)
)
\varphi_{tj\bar{k}}\varphi_{k\bar{p}\bar{l}}}_{B2a} 
\\
&+&
%\left(
(
\underbrace{\varphi_{tj\bar{p}\bar{k}} +C^{t\bar{p}}_{kj}\varphi_{k\bar{k}} 
%\right)
)
%\left(
(
\varphi_{pk\bar{t}\bar{j}} +C^{p\bar{t}}_{jk} \varphi_{k\bar{k}}
%\right)
)}_{B2}
\\
&+&
%\left(
(
\underbrace{\varphi_{pk\bar{t}\bar{j}\bar{p}}}_{\overline{A1}} \,
+ \, C^{p\bar{t}}_{jr}\varphi_{k\bar{r}\bar{p}}  + 
\underbrace{\overline{C^{tp\bar{p}}_{j\bar{k}}}}_{\overline{C1}} 
\varphi_{k\bar{k}}
%\right)
)
\varphi_{tj\bar{k}}
%}_{Hey}
\Big\}  
\end{eqnarray*}
We explain now the meaning of the underbraces. Set $a_{ptjk}:=a_pa_ta_ja_k$ and
\begin{eqnarray*}
B1:=2^5a_{ptjk}Big\{\Big[\varphi_{ptj\bar{k}}-2\sum_la_l
\left(\varphi_{pt\bar{l}}\varphi_{lj\bar{k}}+\varphi_{pl\bar{k}}\varphi_{tj\bar{l}}\right)\Big]\times[\mbox{conjugate}]\Big\}\geq 0
\\
B2:=2^5a_{ptjk}\Big\{\Big(\varphi_{tj\bar{p}\bar{k}}+C^{t\bar{p}}_{kj}\varphi_{k\bar{k}}-2\sum_la_l
\varphi_{tj\bar{l}}\varphi_{l\bar{p}\bar{k}} \Big)\times(\mbox{conjugate})\Big\}\geq 0\,.\qquad
\end{eqnarray*}
Then the underbraced terms in the previous expression of the Laplacian corresponds to the terms $A\ast, C\ast$ of the expression \eqref{CovDerRic} of the covariant derivative of the Ricci tensor and the terms $B\ast$ just defined. To be more precise to see those correspondences we need to make the following change of indexes of the underbraced terms of the Laplacian.
\begin{center}
\begin{tabular}{llllll}
A2& $(k , l , p , t ) \, \rightarrow \, ( l, k , t , p )$ & & \\
A3 & $( k , l , p , t , j ) \, \rightarrow \, ( p , k , j , l , t
)$ & &  
$\overline{\rm A3}$ & $( l , p ) \, \rightarrow \, ( p , l )$
\\
A4 & $( l , j ) \, \rightarrow \, ( j , l )$ 
& & $\overline{\rm A4}$ &
$(t, j, k, l, p ) \, \rightarrow \, ( k, p, l, j, t, )$ \\
A5 & $( k, l, p, j, t, r ) \, \rightarrow \, ( l , k , j , t, r , p)$
&& $\overline{\rm A5}$ & $( l, r ) \, \rightarrow \, )$ \\
A6 & $( k, j, r, t, l ) \, \rightarrow \, ( l, r, k, j, t )$ & &
$\overline{\rm A6}$ & $(r, l, t, j ) \, \rightarrow \, ( l, r, j, t )$
\\
B1a & $( t, j, p, l ) \, \rightarrow \, ( p, l, t, j )$ && B1b & 
$( t, l ) \, \rightarrow \, ( l , t )$ \\
B1c & $( j, l, t, p ) \, \rightarrow \, ( l, j, p, t )$ & &B1d & 
$(t, r, l ) \, \rightarrow \, ( l, t, r )$ 
\\
B2a & $( l , k ) \, \rightarrow \,  ( k ,  l )$ 
\end{tabular}
\end{center} 
Using the expression of the Ricci tensor at the point $x$
$$
\Ric(\omega_{\varphi})_{j\bar{k}}=-ia_p(C^{j\bar{k}}_{pp} +2\varphi_{jp\bar{k}\bar{p}})
+4ia_pa_t\varphi_{jp\bar{t}} \varphi_{t\bar{p}\bar{k}}\,,  
$$
and the expression \eqref{CovDerRic} of the covariant derivative of the Ricci tensor at the point $x$ we find the expression
\begin{eqnarray}
\Delta_{\varphi}|\nabla^{1,0}_{\omega}\partial\bar{\partial}\varphi|^2_{\varphi }
=
-2^5ia_{ltjk}
\Big[
\Ric(\omega_{\varphi})_{l\bar{t}}\,\varphi_{k\bar{j}\bar{l}}\nonumber
+
\Ric(\omega_{\varphi})_{l\bar{j}} \, \varphi_{k\bar{t}\bar{l}} 
+
\Ric(\omega_{\varphi})_{k\bar{l}} \,\varphi_{l\bar{t}\bar{j}}  
\Big]\varphi_{tj\bar{k}}\nonumber
\end{eqnarray}
\begin{eqnarray}\label{C3Expres}
&+&
4\Re e\Big[i\left<\nabla^{1,0}_{\omega}\Ric(\omega_{\varphi}),\nabla^{1,0}_{\omega}\partial\bar{\partial}\varphi\right>_{\varphi}\nonumber
+
2i\left<\Tr_{\varphi}\nabla^{1,0}_{\omega}R_{\omega},\nabla^{1,0}_{\omega}\partial\bar{\partial}\varphi\right>_{\varphi,\omega}\Big]\nonumber
\\
&+&
2^5a_{ptjk}\,
2\Re e\Big[
a_l
\Big(
\underbrace{
C^{j\bar{k}}_{lp}\,\varphi_{tl\bar{p}}
}_{C2}
-
\underbrace{
C^{t\bar{k}}_{pl}\,\varphi_{jp\bar{l}}
}_{C3}
\Big)
\varphi_{k\bar{t}\bar{j}} -4a_lC^{t\bar{p}}_{kj}\,\varphi_{k\bar{k}}  \,
\varphi_{pl\bar{t}}\,\varphi_{k\bar{l}\bar{j}} \nonumber
\\
&+&
C^{t\bar{p}}_{rj} \,\varphi_{rp\bar{k}}\, 
\varphi_{k\bar{t}\bar{j}}  \Big]+B1+B2\,,  \qquad
\end{eqnarray}
where $\Re e$ is the real part of a complex number and
$$
(\Tr_{\varphi}\nabla^{1,0}_{\omega}R_{\omega})(\xi,\eta,\mu):=
\Tr_{\varphi}[\nabla^{1,0}_{\omega,\xi}R_{\omega}(\eta,\mu,\cdot,\cdot)]=
\Tr_{\varphi}[\nabla^{1,0}_{\omega,\xi}R_{\omega}(\cdot,\cdot,\eta,\mu)]\,,
$$
for all $\xi,\eta\in T^{1,0}_{_{X,J}}, \,\mu\in T^{0,1}_{_{X,J}}$.
From now on we reconsider our original notations $\varphi_t=\varphi$ and $\omega_t=\omega_{\varphi}$.
Using the fact that the inverse matrix $(\omega ^{k,\bar{l}}_t)_{k,l}$ evolves by the formula
$$
\frac{d}{d t}\omega_t^{k,\bar{l}}=-\omega_t^{k,\bar{l}}+2\omega_t^{k,\bar{j}}R_{j,\bar{p}}(t) \,\omega_t^{p,\bar{l}},
$$
where $\Ric(\omega_t)=iR_{j,\bar{p}}(t)\,dz_j\wedge d\bar{z}_p$, we find at the point $x$ the expression
\begin{eqnarray}\label{TimeDerNor}
\frac{\partial }{\partial t}|\nabla^{1,0}_{\omega}\partial\bar{\partial}\varphi_t|^2_t
&=&
-3|\nabla^{1,0}_{\omega}\partial\bar{\partial}\varphi_t|^2_t+2^4a_{pljk}
\Big[R_{l\bar{p}}(t)\,\varphi_{k\bar{j}\bar{l}}  +R_{l\bar{j}}(t)\,
\varphi_{k\bar{l}\bar{p}}\nonumber
\\
&+&
R_{k\bar{l}}(t)\,\varphi_{l\bar{p}\bar{j}} \Big]\varphi_{pj\bar{k}} 
+2\Re e\left<\nabla^{1,0}_{\omega}\partial\bar{\partial}\dot{\varphi}_t,\nabla^{1,0}_{\omega}\partial\bar{\partial}\varphi_t\right>_t.
\end{eqnarray}
%%%%%%%%%%%%%%%%%%%%%%%%%%%%%%%%%%%%%%%%%%%%%%%%%%%%%%%%%%%%%%%%%%%%%%%%%%%%%%%%%%%%%%%%%%%%%%%%%%%%%%%%%%%%%%%%%%%%%%%%%%%%%%%%%%%%%%%%%%%
% \omega^{\bar{}}_{\varphi}  _{}  \varphi_{\bar{}}  \varphi_{ \bar{}\bar{}}  
%                                                                     \varphi_{ \bar{}\bar{}\bar{}}
%%%%%%%%%%%%%%%%%%%%%%%%%%%%%%%%%%%%%%%%%%
Using the expression \eqref{C3Expres} and the fact that all the metrics $\omega_t$ are uniformly equivalents to the initial metric $\omega$,  we obtain the inequality
\begin{eqnarray}\label{C3LapKRcci}
\Delta _t|\nabla^{1,0}_{\omega}\partial\bar{\partial}\varphi_t|^2_t
\geq
2^5a_{ptjk}
\Big[R_{l\bar{p}}(t)\,\varphi_{k\bar{j}\bar{l}}  +R_{l\bar{j}}(t)\,
\varphi_{k\bar{l}\bar{p}}
+
R_{k\bar{l}}(t)\,\varphi_{l\bar{p}\bar{j}} \Big]\varphi_{pj\bar{k}} \nonumber
\\
+
4\Re e\left<i\nabla^{1,0}_{\omega}\Ric(\omega_t),
\nabla^{1,0}_{\omega}\partial\bar{\partial}\varphi_t\right>_t-C_1
|\nabla^{1,0}_{\omega}\partial\bar{\partial}\varphi_t|^2_t 
-
C'_2
|\nabla^{1,0}_{\omega}\partial\bar{\partial}\varphi_t|_t,
\end{eqnarray}
where $C_1,\,C'_2>0$ are two constants independents of $t$. By deriving the K\"ahler-Ricci flow identity $i\partial\bar{\partial}\dot{\varphi}_t=\omega_t-\Ric(\omega_t)$, we find the equality 
$$
i\nabla^{1,0}_{\omega}\partial\bar{\partial}\dot{\varphi}_t=i\nabla^{1,0}_{\omega}\partial\bar{\partial}\varphi_t-\nabla^{1,0}_{\omega}\Ric(\omega_t),
$$
which combined with the relations \eqref{TimeDerNor} and \eqref{C3LapKRcci} gives the uniform estimate
\begin{eqnarray*}
\Openbox_{\,t}|\nabla^{1,0}_{\omega}\partial\bar{\partial}\varphi_t|^2_t
\geq (2-C_1)|\nabla^{1,0}_{\omega}\partial\bar{\partial}\varphi_t|^2_t-C'_2
|\nabla^{1,0}_{\omega}\partial\bar{\partial}\varphi_t|_t
\geq
-C_1
|\nabla^{1,0}_{\omega}\partial\bar{\partial}\varphi_t|^2_t -C_2\,,
\end{eqnarray*}
for some uniform constant $C_2>0$ sufficiently big. \hfill$\Box$
%%%%%%%%%%%%%%%%%%%%%%%%%%%%%%%%%%%%%%%%%%%%%%%%%%%%%%%%%%%%%%%%%%%%%%%%%%%%%%%%%%%%%%%%%%%%%%%%%%%%%%%%%%%%%%%%%%%%%%%%%%%%%%%%%%%%%%%%%%%%%%%%%%%%%%%%%%%%%%%%%%%%%%%%%%%%%%%%%%%%%%%%%%%%%%%%%%%%%%%%%%%%%%%%%%%%%%%%%%%%%%%%%%%%%%
\section{The existence of a K\"ahler-Einstein metric}
%%%%%%%%%%%%%%%%%%%%%%%%%%%%%%%%%%%%%%%%%%%%%%%%%%%%%%%%%%%%%%%%%%%%%%%%%%%%%%%%%%%%%%%%%%%%%%%%%%%%%%%%%%%%%%%%%%%%%%%%%%%%%%%%%%%%%%%%%%%%%%%%%
We remind that the uniform estimate $\omega_t^n\geq k\,\omega^n$ is equivalent to the uniform estimate $\varphi_t+c_t\leq C$. Then the identity \ref{KenKRfloNorm} implies that the K-energy is also uniformly bounded from below along the flow, thus the limit $\lim_{t\rightarrow +\infty}\nu_{\omega}(\varphi_t)$ is finite.
We remind also that along the K\"ahler-Ricci flow we have the identity
$$
\frac{d}{d t}\nu_{\omega}(\varphi_t)=-\dashint\limits_X|\partial \dot{\varphi}_t|^2_t\,\omega^n_t\,.
$$
%which implies the convergence of the integral$$\int\limits_0^{+\infty}dt\;\dashint\limits_X|\partial \dot{\varphi}_t|^2_t\,\omega^n_t=\nu_{\omega}(\varphi_0)-\lim_{T\rightarrow +\infty}\nu_{\omega}(\varphi_T)\leq \nu_{\omega}(\varphi_0)+C<+\infty\,.$$
So for all increasing sequences of times 
$(\tau_k)\subset [0,+\infty),\;\tau_k\rightarrow +\infty$ there exist a sequence $(t_k)\,,t_k\in [\tau_k,\tau_{k+1}]$ such that
\begin{eqnarray}\label{Enrg0seq}
\lim_{k\rightarrow +\infty}\,\dashint\limits_X|\partial \dot{\varphi}_{t_k}|^2_{t_k}\,\omega^n_{t_k}=0.
\end{eqnarray}
Moreover the $C^2$ and $C^3$-uniform estimates $|\partial\bar\partial \varphi_t|_{C^0(X)},\,|\nabla^{1,0}_{\omega}\partial\bar\partial \varphi_t|_{C^0(X)}\leq C$ implient that the $(1,1)$-formes $(\partial\bar\partial \varphi_t)_{t\in [0,+\infty)}$ are uniformly bounded in the $C^{\alpha}(X)$-topology. 
The operator $\Delta_t$ is uniformly elliptic with coefficients uniformly bounded in $C^{\alpha}$-norm, at least. The right hand side of the equation \eqref{C3Alfa} 
\begin{eqnarray*}
\Openbox_{\,t}(\xi.\varphi_t)+2\xi .\varphi_t=(\Tr_{\omega}-\Tr_t)(L_{\xi}\,\omega)+2\xi . h_{\omega}\,,
\end{eqnarray*}
$\xi\in{\cal E}(T_X)(U)$,
is also uniformly bounded in $C^{\alpha}$-norm, at least. By the regularity theory for parabolic equations \cite{Lad} we deduce that the functions $(\xi .\varphi_t)_{t\in [0,+\infty)}$ are uniformly bounded in $C^{2,\alpha}$-norm. Then the $C^0$-uniform estimate $|\varphi_t+c_t|\leq C$ implies the existence of a subsequence $(s_k)$ of $(t_k)$ such that the sequences $(\varphi_{s_k}+c_{s_k}),\,(d\varphi_{s_k}),\,(\partial\bar\partial \varphi_{s_k})$ and $(\nabla^{1,0}_{\omega}\partial\bar\partial\varphi_{s_k})$ convergent uniformly respectively to $\varphi_{\infty},\,d\varphi_{\infty},\,\partial\bar\partial\varphi_{\infty},\,$ and $\nabla^{1,0}_{\omega}\partial\bar\partial\varphi_{\infty}$. The uniform estimate $\omega^n_t/\omega^n\geq K_0>0$ gives $\omega^n_{\varphi_{\infty}}/\omega^n\geq K_0>0$, which implies $i\partial\bar\partial\varphi_{\infty}>-\omega$. Moreover we deduce the existence of the limits
$$
\psi:=\lim_{k\rightarrow +\infty}\dot{\varphi}_{s_k}=\log\frac{\omega^n_{\varphi_{\infty}}}{\omega^n}+
\varphi_{\infty}-h_{\omega}
$$
and $\partial \psi=\lim_{k\rightarrow +\infty}\partial\dot{\varphi}_{s_k}$ in the topology of the uniform convergence at least. Then the limit \eqref{Enrg0seq} implies 
$$
0=\lim_{k\rightarrow +\infty}\,\dashint\limits_X|\partial \dot{\varphi}_{s_k}|^2_{s_k}\,\omega^n_{s_k}=\dashint\limits_X|\partial \psi|^2_{\varphi_{\infty}}\omega^n_{\varphi_{\infty}},
$$
which means $\psi=0$, by the integral normalization of $\dot{\varphi}_t$. So we have a solution $\varphi_{\infty}\in {\cal P}_{\omega}^{3,\alpha}$ of the elliptic non-linear equation
$$
F(\varphi_{\infty}):=\log\frac{\omega^n_{\varphi_{\infty}}}{\omega^n}+
\varphi_{\infty}-h_{\omega}=0.
$$
The ellipticity follows from the fact that $i\partial\bar\partial\varphi_{\infty}>-\omega$ and the expression of the differential 
$d_{\varphi_{\infty}}F(v)=2^{-1}\Delta_{\varphi_{\infty}}v+v$. By Schauder elliptic regularity (see \cite{Aub}, Th. 3.56, pag. 86) we deduce that the solution $\varphi_{\infty}$ is smooth. In conclusion we have solve the Einstein equation $\Ric(\omega_{\varphi_{\infty}})=\omega_{\varphi_{\infty}}$.
Clearly the Einstein metric $\omega_{\varphi_{\infty}}$ is $G$-invariant if the initial metric $\omega$ of the K\"ahler-Ricci flow is $G$-invariant.\hfill$\Box$
%%%%%%%%%%%%%%%%%%%%%%%%%%%%%%%%%%%%%%%%%%%%%%%%%%%%%%%%%%%%%%%%%%%%%%%%%%%%%%%%%%%%%%%%%%%%%%%%%%%%%%%%%%%%%%%%%%%%%%%%%%%%%%%%%%%%%%%%%%%%%%%%%%%%%%%%%%%%%%%%%%%%%%%%%%%%%%%%%%%%%%%%%%%%%%%%%%%%%%%%%%%%%%%%%%%%%%%%

\vspace{1cm}
\noindent
Nefton Pali
\\
Mathematics Department - Princeton University
\\
Fine Hall - Washington road
\\
08544 Princeton, NJ - USA
\\
E-mail: \textit{npali@math.princeton.edu}
\end{document}